\documentclass[a4paper,11pt]{article}
\usepackage{makeidx}
\usepackage{graphicx}
\usepackage{latexsym}
\usepackage{amsmath}
\usepackage{amscd}
\usepackage{amssymb}

\usepackage{a4wide}
\author{St\'ephane Launois \footnote{This work was partially supported by Leverhulme Research Interchange Grant F/00158/X}
 \\
\\
{\small{\it Laboratoire de Math\'ematiques - UMR6056, Universit\'e de Reims}}\\
\small{{\it Moulin de la Housse - BP 1039 - 51687 REIMS Cedex 2, France}}\\
\small{e-mail : stephane.launois@univ-reims.fr}}

\title{Combinatorics of $\hc$-primes in quantum matrices.}

 \date{ }

\newcommand{\fin}{$\blacksquare$}
\newcommand{\preuve}{\textit{Proof.} }
\newcommand{\gc}{ [ \hspace{-0.65mm} [}
\newcommand{\dc}{]  \hspace{-0.65mm} ]}
\newcommand{\ov}{\overline}

\newcommand{\mnc}{O_{q}\left( \mathcal{M}_n(\mathbb{C}) \right)}
\newcommand{\mmc}{O_{q}\left( \mathcal{M}_m(\mathbb{C}) \right)}
\newcommand{\mmpc}{O_{q}\left( \mathcal{M}_{m,p}(\mathbb{C}) \right)}

\newcommand{\sn}{O_{q}\left( SL_n(\mathbb{C}) \right) }

\newcommand{\gn}{O_{q}\left( GL_n(\mathbb{C}) \right) }

\newcommand{\ia}{i,\alpha}

\newcommand{\spec}{\mathrm{Spec}}
\newcommand{\fract}{\mathrm{Fract}}

\newcommand{\yia}{Y_{i,\alpha}}

\newcommand{\hc}{\mathcal{H}}

\newcommand{\xia}{X_{i,\alpha}}

\begin{document}
\maketitle

\newtheorem{theo}{Theorem}[section]
\newtheorem{defi}[theo]{Definition}

\newtheorem{rem}[theo]{Remark}

\newtheorem{nota}[theo]{Notation}

\newtheorem{prop}[theo]{Proposition}
\newtheorem{hypo}[theo]{Hypothesis}
\newtheorem{lem}[theo]{Lemma}
\newtheorem{conv}[theo]{Convention}

\newtheorem{cor}[theo]{Corollary}
\newtheorem{obs}[theo]{Observation}

\begin{abstract}
$ $

For $q \in \mathbb{C}$ generic, we give an algorithmic construction of an ordered bijection between the set of $\hc$-primes of $\mnc$ and the sub-poset $\mathcal{S}$ of the (reverse) Bruhat order of the symmetric group $S_{2n}$ consisting of those permutations that move any integer by no more than $n$ positions. Further, we describe the permutations that correspond via this bijection to rank $t$ $\hc$-primes. 
More precisely, we establish the following result. Imagine that there is a barrier between positions $n$ and $n+1$. Then a $2n$-permutation $\sigma  \in \mathcal{S}$ corresponds to a rank $t$ $\hc$-invariant prime ideal of $\mnc$ if and only if the number of integers that are moved by $\sigma$ from the right to the left of this barrier is exactly $n-t$. The existence of such a bijection was conjectured by Goodearl and Lenagan.
\end{abstract}
$ $
\\$ $
\\$ $

\begin{center} \begin{large}\textbf{Introduction}\end{large} \end{center}
$ $

Let $m$ and $p$ be two integers greater than or equal to $2$, and let $q$ be 
a complex number transcendental over $\mathbb{Q}$. Denote 
by $A=\mmpc$ the quantization of the ring of regular functions on 
$m \times p$ matrices with entries in $\mathbb{C}$ and by $(Z_{\ia})_{(\ia) \in \gc 1,m \dc \times \gc 1,p \dc}$ 
the matrix of its canonical generators. It is well known that the group 
$\hc:=(\mathbb{C}^*)^{m+p}$ acts on $A$ by $\mathbb{C}$-automorphisms via:
$$(a_1,\dots,a_m,b_1,\dots,b_p).Z_{\ia}=a_i b_{\alpha} Z_{\ia} 
\quad ((\ia) \in \gc 1,m \dc  \times \gc 1,p \dc).$$

In \cite{gl1}, Goodearl and Letzter have shown that the set $\hc$-$\spec(A)$ of all $\hc$-invariant prime ideals of $A$ is finite and that, in order to calculate the prime and primitive spectra of $R$, it is enough to determine the $\hc$-invariant prime ideals of $A$. More precisely, they have shown that the prime spectrum of $A$ admits a natural stratification, indexed by the finite set $\hc$-$\spec(A)$, such that each stratum is Zariski-homeomorphic to the prime spectrum of a commutative Laurent polynomial ring over $\mathbb{C}$. Further, let us mention that this stratification is a powerful tool to recognize the primitive ideals since Goodearl and Letzter have also proved that the primitive ideals of $A$ are precisely the primes of $A$ that are maximal in their stata (see \cite{gl1}).

Hence, to understand the prime and primitive spectra of $A$, the first current step is to study these 
$\hc$-invariant prime ideals. Several progress have been made recently to understand these ideals. First, in the case where $m=p$, Goodearl and Lenagan have shown 
that, in order to obtain descriptions of all the $\hc$-invariant prime ideals of $A$, we just need to 
determine the $\hc$-invariant prime ideals of certain "localized step-triangular factors" of $A$ (see \cite[Theorem 3.5]{glen1}). Using this result, Goodearl and Lenagan have 
computed the $\hc$-invariant prime ideals of $O_{q}\left( \mathcal{M}_2(\mathbb{C}) \right)$ 
(see \cite{glen1}) and $O_{q}\left( \mathcal{M}_3(\mathbb{C}) \right)$ (see \cite{glen2}). 
Next, using the theory of deleting-derivations (see \cite{c1}), Cauchon has given a description of the set $\hc$-$\spec(A)$. As a consequence of this description, he has obtained a formula for the total number of 
$\hc$-invariant prime ideals in $A$ (see \cite[Proprosition 3.3.2]{c2}). Still using Cauchon's theory of deleting-derivations, we have proved (see \cite{lau}) that the $\hc$-invariant prime ideals in $A$ are generated by quantum minors, as conjectured by Goodearl and Lenagan. This result has allowed us to construct 
an algorithm which provides an explicit generating set of quantum minors for each $\hc$-invariant prime ideal in $A$ (see \cite{lau2}).

In this paper we investigate the combinatoric of the set $\hc$-$\spec(A)$. More precisely, Goodearl and Lenagan have conjectured that, in the case where $m=p$, there exists an explicit order-preserving bijection between $\hc$-$\spec(\mmc)$ and the sub-poset of the (reverse) Bruhat order of $S_{2m}$ consisting of those permutations that move any integer by no more than $m$ positions. The main goal of this paper is to construct such a bijection. In fact, we will construct an order-preserving bijection between $\hc$-$\spec(A)$ and the sub-poset 
$\mathcal{S}:=\{ \sigma \in S_{m+p} \mid -p \leq i - \sigma(i) \leq m \mbox{ for all } i \in \gc 1,m+p \dc \}$ of the (reverse) Bruhat order of $S_{m+p}$. To do this, we will use the $(m,p)$ deleting-derivations algorithm that we have introduced in \cite{lau} to prove that the $\hc$-invariant prime ideals of $A$ are generated by quantum minors. Set $n:=m+p$. This algorithm consists of certain changes of variables in the field of fractions of $R:=\mnc$. It is interesting to note that, at the step $(m,m)$ of this algorithm, the subalgebra $R^{(m,m)}$ of $\fract(R)$ generated by the new indeterminates is isomorphic to an iterated Ore extension of $A$ that does not involve $\sigma$-derivations. This fact allows us to define an embedding from $\hc$-$\spec(A)$ into a subset of $\spec(R^{(m,m)})$, the image of $J \in \hc$-$\spec(A)$ just being the ideal of $R^{(m,m)}$ generated by $J$. Next, using some suitable localizations and contractions related to the involved changes of variables, we are able to extend this embedding to an embedding from $\hc$-$\spec(A)$ into the set of those primes of $R$ that are invariant under the natural action of $\hc_R:=(\mathbb{C}^*)^{2n}$ on $R$ and that do not contain the quantum determinant of $R$. Note that the construction of this embedding is in fact algorithmic and is done with a finite number of simple steps. Finally, using classical results of localization theory together with the isomorphism between 
$\sn[z,z^{-1}] $ and $\gn$ that was constructed by Levasseur and Stafford (see \cite{ls}), we obtain an order-preserving embedding $\psi$ from $\hc$-$\spec(A)$ into the set $\hc'_R$-$\spec(\sn)$ of those primes of the quantum group $\sn$ that are invariant under the natural action of the group 
$ \hc'_R:=\{(a_1,\dots,a_n,b_1,\dots,b_n) \in (\mathbb{C}^*)^{2n} \mid a_1  \dots a_n b_1 \dots b_n =1\}$ on $\sn$.

Since $q$ is transcendental over $\mathbb{Q}$, the set $\hc'_R$-$\spec(\sn)$ is well known. Using the partition of $\spec(\sn)$ constructed by Hodges and Levasseur (see \cite{hl1,hl2}), and a theorem of Joseph (see \cite{jos2}), Brown and Goodearl have shown (see \cite{bg}) that this set is parametrized by $S_n \times S_n$, and that, for each $w \in  S_n \times S_n$, the corresponding ideal $I_w$ is generated by some explicit quantum minors. So, in order to compute the image of $\psi$, we first establish some technical criteria for a quantum minor to belong to a prime ideal in the image of $\psi$. Then these tools allow us to prove that the image of $\psi$ 
is contained in the set of those $I_{w_0,w_0\sigma}$ such that $\sigma \in  \mathcal{S}=\{ \sigma \in S_{m+p} \mid -p \leq i - \sigma(i) \leq m \mbox{ for all } i \in \gc 1,n \dc \}$ (where $w_0$ denotes the longest element of $S_n$). Finally, we observe that these two sets have the same cardinality, so that they are actually equal. As a consequence, $\psi$ induces an ordered bijection between $\hc$-$\spec(A)$ and the set of those $I_{w_0,w_0\sigma}$ with $\sigma \in  \mathcal{S}$. Naturally the inverse of this bijection, which is also an order-preserving bijection, allows us to construct an ordered bijection $\Xi$ from the sub-poset $\mathcal{S}$ of the (reverse) Bruhat order of $S_n$ and $\hc$-$\spec(A)$.

Note that, very recently, Brown, Goodearl and Yakimov (see \cite{bgy}) have investigated the geometry of the Poisson variety $\mathcal{M}_{m,p}(\mathbb{C})$, where the Poisson structure is the standard one that can be obtained from commutators of $\mmpc$. (This Poisson structure must be viewed as the semiclassical limit of the quantization $\mmpc$.) Their results together with the bijection $\Xi$ constructed in this paper prove the existence of a bijection between $\hc$-$\spec(\mmpc)$ and the set of $\hc$-orbits of symplectic leaves in $\mathcal{M}_{m,p}(\mathbb{C})$. 

In the last section of the paper, assuming that $m=p$, we describe the permutations that correspond via this bijection $\Xi$ to rank $t$ $\hc$-primes in $A$, that is, to those $\hc$-invariant prime ideals of $A$ which contain all $(t+1) \times (t+1)$ quantum minors but not all $t \times t$ quantum minors. More precisely, we establish the following result. Imagine that there is a barrier between positions $m$ and $m+1$. Then a $2m$-permutation $\sigma$ belonging to $\mathcal{S}$ corresponds to a rank $t$ $\hc$-invariant prime ideal of $A$ if and only if the number of integers that are moved by $\sigma$ from the right to the left of this barrier is exactly $m-t$. The existence of an order-preserving bijection from the sub-poset $\mathcal{S}$ of the (reverse) Bruhat order of $S_{2m}$ and $\hc$-$\spec(A)$ with such properties was conjectured by Goodearl and Lenagan.

Note that the question of such a construction remains open for a general base field $\mathbb{K}$ and a non-zero $q \in \mathbb{K}$ not a root of unity. However we have to say that, if we restrict to $\mathbb{K} = \mathbb{C}$ and $q$ transcendental over $\mathbb{Q}$, it is only for a technical condition: the ideals $I_w$ ($w \in S_n \times S_n$) introduced by Hodges and Levasseur in \cite{hl1} are known to be prime only in this case (see \cite{jos2}).
\\$ $

Throughout this paper, we use the following conventions.
\\$ $
\\$\bullet$ If $I$ is a finite set, $|I|$ denotes its cardinality.
\\$\bullet$ $\mathbb{C}$ denotes the field of complex numbers and we set $\mathbb{C}^*:=\mathbb{C}\setminus \{0\}$.
\\$\bullet$ $q\in \mathbb{C}^*$ is transcendental over $\mathbb{Q}$.
\\$\bullet$ $m,p$ denote two positive integers with $m,p \geq 2$, and we set $n:=m+p$.
\\$\bullet$ $A=\mmpc$ denotes the quantization of the ring of regular functions on 
$m \times p$ matrices with entries in $\mathbb{C}$; it is the $\mathbb{C}$-algebra generated by the 
$m \times p$ indeterminates $Z_{\ia}$, $1 \leq i \leq m$ and $ 1 \leq \alpha \leq p$, subject to the following 
relations.\\
If $\left( \begin{array}{cc} x & y \\ z & t \end{array} \right)$ is any $2 \times 2$ sub-matrix of 
$\mathcal{Z}:=\left( Z_{\ia} \right) _{(\ia) \in \gc 1,m \dc \times \gc 1,p \dc}$, then
 \begin{enumerate}
 \item $ yx=q^{-1} xy, \quad zx=q^{-1} xz, \quad zy=yz, \quad ty=q^{-1}
 yt, \quad tz=q^{-1} zt$.
 \item $tx=xt-(q-q^{-1})yz$.
\end{enumerate}
 It is well known that $A$ can be presented as an iterated Ore extension over $\mathbb{C}$, with the generators $Z_{\ia}$ adjoined in lexicographic order. Thus the ring $A$ is a Noetherian domain. Moreover, since 
$q$ is transcendental over $\mathbb{Q}$, it follows from \cite[Theorem 3.2]{gl4} that all prime ideals 
of $A$ are completely prime.
\\$\bullet$ It is well known that the group $\hc:=\left( \mathbb{C}^* \right)^{m+p}$ acts 
on $A$ by $\mathbb{C}$-algebra automorphisms via:
$$(a_1,\dots,a_m,b_1,\dots,b_p).Z_{\ia} = a_i b_\alpha Z_{\ia} 
\quad \forall \: (\ia)\in \gc 1,m \dc \times \gc 1,p \dc.$$
An ideal $I$ of $A$ is said to be \textbf{$\hc$-invariant} if $h.I =I$ for all $h\in \hc$. 
We denote by $\hc$-$\spec(A)$ the set of $\hc$-invariant prime ideals of $A$. Recall (see \cite[5.7. (i)]{gl1}) that $A$ has only finitely many $\hc$-invariant prime ideals.\\$ $


\section{Cardinality of some subsets of $S_n$.}
\label{sectionbruhatS}
$ $

Recall that $n=m+p$. We denote by $S_n$ the group of permutations of $\gc 1,n \dc$ and by $w_0$ the longest element of $S_n$. Recall that $w_0(i) = n+1 -i$ for all $i \in \gc 1,n \dc$ and that $w_0=s_{1,n} s_{2,n-1} \dots$, where $s_{i,j}$ denotes the transposition associated to $i$ and $j$ ($i \neq j$).

In this section, we are going to study the following sub-poset of the (reverse) Bruhat ordering of $S_n$:
$$\mathcal{S} := \{\sigma \in S_n \mid \sigma \leq \sigma_0\},$$
where $\sigma_0$ is the $n$-permutation defined by 
$$ \sigma_0 (i) = \left\{ \begin{array}{ll} p+i & \mbox{ if } i \in \gc 1,m \dc \\
i-m & \mbox{ else.} 
\end{array} \right. $$
In fact, our main aim in this section is to prove that the cardinality of $\mathcal{S}$ is equal to the poly-Bernoulli number $B_p^{(-m)}$ (we send back to \cite{kaneko1} for a precise definition of poly-Bernoulli numbers). To do this, we will first rewrite the sub-poset $\mathcal{S}$ of the (reverse) Bruhat order of $S_n$ as a subset of $S_n$ of restricted permutations. More precisely, we will show that 
$$\mathcal{S} = \{ \sigma \in S_n \mid -p \leq i - \sigma(i) \leq m \mbox{ for all } i \in \gc 1,n \dc \}.$$
Next we will rewrite the formula obtained by Vesztergombi for the cardinality of this set of restricted permutations (see \cite{ves}) in order to show that $\mid \mathcal{S} \mid =B_p^{(-m)}$. As a corollary, we obtain that $\mathcal{S}$ has the same cardinality as $\hc$-$\spec(A)$. This result will be a key-point in order to construct in the sequel an explicit order-preserving bijection from $\mathcal{S}$ onto $\hc$-$\spec(A)$.

Finally, in the case where $m=p$, we also compute the cardinality of some subsets $\mathcal{S}_t$ ($t \in \gc 0, m \dc$) of $\mathcal{S}$. We will see in the last section of this paper (see \ref{sectionrangtperm}) that these subsets of $\mathcal{S}$ are closely related to the rank $t$ $\hc$-primes of $A$.

\subsection{The Bruhat ordering on $S_n$.}
$ $

To avoid any ambiguity, we need to make precise our convention for the Bruhat order. In this paper, we will always work with the reverse of the usual Bruhat ordering introduced in \cite{Dixmier}, so that $id \leq \sigma \leq w_0$ for all $\sigma \in S_n$. 

\begin{nota}
Let $\sigma,\sigma' \in S_n$ and $j \in \gc 1,n-1 \dc$.
\\If $\sigma(\gc 1,j \dc)=\{\sigma_1,\dots,\sigma_j\}$ with $\sigma_1 < \dots < \sigma_j$ and $\sigma'(\gc 1,j \dc)=\{\sigma'_1,\dots,\sigma'_j\}$ 
with $\sigma'_1 < \dots < \sigma'_j$, then we write $\sigma \leq_j \sigma'$ if $\sigma_k \leq \sigma'_k$ for all $k \in \gc 1, j \dc$.
\end{nota}

Recall (see \cite[Proposition 1.8.1]{hl1}) that the (reverse) Bruhat ordering can be characterized as follows. 

\begin{prop}
\label{propbruhat}
Let $\sigma,\sigma' \in S_n$. Then $\sigma \leq \sigma'$ if and only if $\sigma \leq_j \sigma'$ for all $j \in \gc 1, n-1 \dc$.
\end{prop}

\subsection{The cardinality of $\mathcal{S}$.}
$ $

Recall that $n=m+p$. Before computing the cardinality of $\mathcal{S}$, we give another description of this set in terms of restricted permutations. Note that this result also appears in \cite[Lemma 3.12]{bgy}.

\begin{prop}
\label{descriptionS}
$\mathcal{S} = \{ \sigma \in S_n \mid -p \leq i - \sigma(i) \leq m \mbox{ for all } i \in \gc 1,n \dc \}$.
\end{prop}
\preuve First, let $\sigma \in S_n$ such that $\sigma \leq \sigma_0$. Then it follows from Proposition \ref{propbruhat} that, for all $j \in \gc 1, n-1 \dc$, 
we have $\sigma \leq_j \sigma_0$. We have to prove that, for all $i \in \gc 1,n \dc$, 
we have $-p \leq i - \sigma(i) \leq m$. To do this, we distinguish two cases.

$\bullet$ Assume that $i\in \gc 1,m \dc$. Then it is clear that we have $i - \sigma(i) \leq m$. On the other hand, since $\sigma \leq_i \sigma_0$ and $\sigma_0(\gc 1,i \dc)=\{p+1,\dots,p+i\}$, we have $\sigma(i) \leq p+i$. Thus $-p \leq i-\sigma(i) $ as required. 

$\bullet$ We now assume that $i \in \gc m+1,n \dc$. Then we have $i - \sigma(i) \geq m+1 -n =1-p \geq -p$. On the other hand, since $\sigma \leq_{i-1} \sigma_0$ and $\sigma_0( \gc 1,i-1 \dc ) = \gc 1,i-1-m \dc \cup \gc p+1, n \dc$, there exist 
$k_{1},\dots,k_{i-1-m} \in \gc 1,i-1 \dc$ such that $\sigma(k_{1})=1$, ..., $\sigma(k_{i-1-m})=i-1-m$. 
Hence $\sigma(i) \geq i-m$ and so $i-\sigma(i) \leq m$, as desired.

So we have just proved that $ \mathcal{S} \subseteq  \{ \sigma \in S_n \mid -p \leq i - \sigma(i) \leq m \mbox{ for all } i \in \gc 1,n \dc \}$.

Let now $\sigma \in  \{ \tau \in S_n \mid -p \leq i - \tau(i) \leq m \mbox{ for all } i \in \gc 1,n \dc \}$. We need to show that $\sigma \leq \sigma_0$, that is, in view of Proposition \ref{propbruhat}, $\sigma \leq_i \sigma_0$ for all $i\in\gc 1,n-1 \dc$. Once again, we distinguish two cases.

$\bullet$ Assume that $i \in \gc 1,m \dc$. Since $\sigma \in \{ \tau \in S_n \mid -p \leq k -  \tau(k) \leq m \mbox{ for all } k \in \gc 1,n \dc \}$, we have $ \sigma(j) \leq p+j =\sigma_0(j)$ for all $j \in \gc 1,i \dc$, so that 
$\sigma \leq_i \sigma_0$.

$\bullet$ Assume now that $i \in \gc m+1,n-1 \dc$. Since $\sigma_0 (\gc 1,i \dc) = \gc 1,i-m \dc \cup \gc p+1, n \dc$, to prove that $\sigma \leq_i \sigma_0$, it is sufficient to show that $\gc 1, i-m \dc \subseteq  \sigma  (\gc 1,i \dc)$. Let $j \in \gc i+1,n \dc$. Since $\sigma \in \{ \tau \in S_n \mid -p \leq k - \tau(k) \leq m \mbox{ for all } k \in \gc 1,n \dc \}$, we have $\sigma(j) \geq j-m \geq i+1-m$. Thus $\sigma (\gc i+1,n \dc) \subseteq \gc i+1-m, n \dc$. This implies that $\gc 1, i-m \dc \subseteq \sigma  (\gc 1,i \dc)$ as desired. \fin
\\$ $

We deduce from the above description of the set $\mathcal{S}$ that 
\begin{eqnarray*} 
\mid \mathcal{S} \mid &  = &  \mid \{ \sigma \in S_n \mid -p \leq i - \sigma(i) \leq m \mbox{ for all } i \in \gc 1,n \dc \} \mid .
 \end{eqnarray*}
Now, the cardinality of this last set is known. In \cite[Theorem 4]{ves}, Vesztergombi has computed this number and has obtained the following formula. 
$$ \mid \{ \sigma \in S_n \mid -p \leq i - \sigma(i) \leq m \mbox{ for all } i \in \gc 1,n \dc \} \mid =  2! \sum_{i=0}^{p-1} (-1)^{p-1+i} i! \left( \begin{array}{c} 2+i \\ i \end{array} \right)  (2+i)^{m-1} S(p,i+1),$$ 
where $S(l,k)$ denotes the Stirling number of second kind (see, for example, \cite{stanleybook} for more details on the Stirling numbers of second kind). Hence we have 
$$ \mid   \mathcal{S} \mid =  \sum_{i=0}^{p-1} (-1)^{p-1+i} (i+1)!  (2+i)^{m} S(p,i+1),$$ 
that is, 
$$ \mid   \mathcal{S} \mid = (-1)^p \sum_{i=1}^{p} (-1)^{i} i!  (i+1)^{m} S(p,i).$$
Observing that $S(p,0)=0$, this leads to 
 $$ \mid   \mathcal{S} \mid = (-1)^p \sum_{i=0}^{p} (-1)^{i} i!  (i+1)^{m} S(p,i),$$
and so we deduce from \cite[Theorem 1]{kaneko1} that 
$$ \mid   \mathcal{S} \mid = B_p^{(-m)},$$
where $B_p^{(-m)}$ denotes the poly-Bernoulli number associated to $p$ and $-m$ (see \cite{kaneko1} for the definition of poly-Bernoulli numbers). To resume, we have just proved the following statement.

\begin{prop}
\label{cardinalityS}
The cardinality of $\mathcal{S}$ is given by the poly-Bernoulli number $B_p^{(-m)}$($=B_m^{(-p)}$).
\end{prop}

Recall that $A= \mmpc$. To prove that $\mid \mathcal{S} \mid = \mid \hc \mbox{-}\spec(A) \mid$, it remains to prove that $\mid \hc \mbox{-}\spec(A) \mid$ is also equal to the poly-Bernoulli number $B_p^{(-m)}$. Note that, in the case where $m=p$, Goodearl and McCammond have rewritten the formula obtained by Cauchon for the total number of $\hc$-invariant prime ideals in $A$ in order to prove this result (see \cite[2.2]{lau4}). In the general situation, this can be proved by similar arguments. As in \cite[3.3]{c2}, we show that 
$$\mid \hc \mbox{-}\spec(A) \mid  = (-1)^{p-1} \sum_{k=2}^{p+1} k^m \sum_{j=1}^{k-1} (-1)^{j-1} \left( \begin{array}{c} k-1 \\ j \end{array} \right) j^p .$$
Next, as in \cite[2.2]{lau4}, we get $\mid \hc \mbox{-} \spec(A) \mid = B_p^{(-m)}$, and so we deduce from Proposition \ref{cardinalityS} the following result.

\begin{cor}
\label{cardok}
$\mid \mathcal{S} \mid = \mid \hc \mbox{-}\spec(A) \mid$.
\end{cor} 

\subsection{Cardinality of some subsets of $\mathcal{S}$ when $m=p$.}
\label{sectionST}
$ $

In this paragraph, we assume that $m=p$, so that we have 
$$\mathcal{S}= \{ \sigma \in S_{2m} \mid -m \leq  \sigma(i) -i  \leq m \mbox{ for all } i \in \gc 1,2m \dc \}.$$
Now, imagine that there is a barrier between positions $m$ and $m+1$. The aim of this section is to count the $2m$-permutations $\sigma \in \mathcal{S}$ such that the number of integers that are moved by $\sigma$ from the right to the left of this barrier is exactly $t$ ($t \in \gc 0,m \dc$). 

\begin{nota} 
For all $t \in \gc 0,m \dc$, we denote by $\mathcal{S}_t$ the set of those $\sigma \in \mathcal{S}$ such that there exist $m+1 \leq j_1 < \dots < j_t \leq 2m$ with 
\begin{enumerate}
\item $\sigma(j_k) \leq m \mbox{ for all } k \in \gc 1,t \dc$,
\item $ \mbox{ and } \sigma (k) \geq m+1 
\mbox{ for all } k \in \gc m+1,2m \dc \setminus \{j_1,\dots,j_t\}$.
\end{enumerate}
\end{nota}
$ $

\begin{prop}
\label{cardinalSt}
For all $t \in \gc 0,m \dc$, we have $\mid \mathcal{S}_t \mid =  \left[ (m-t)! S(m+1,m-t+1) \right]^2$.
\end{prop} 
\preuve First, observe that  $\sigma \in \mathcal{S}_t$ if and only if there exist $m+1 \leq j_1 < \dots < j_t \leq 2m$ and $1 \leq i_1 < \dots < i_t \leq m$ with  
\begin{enumerate}
\item $\sigma(j_k) \leq m \mbox{ for all } k \in \gc 1,t \dc$,
\item $ \sigma (l) \geq m+1 
\mbox{ for all } l \in \gc m+1,2m \dc \setminus \{j_1,\dots,j_t\} $,
\item $\sigma(i_k) \geq m+1 \mbox{ for all } k \in \gc 1,t \dc$,
\item $ \sigma (l) \leq m 
\mbox{ for all } l \in \gc 1,m \dc \setminus \{i_1,\dots,i_t\} $
\item and $\mid \sigma(l) -l \mid \leq m$ for all $k \in \gc 1,2m \dc$.
\end{enumerate}
Note that the last condition is empty when $l \in \gc 1, 2m \dc \setminus \{ i_1,\dots,i_t,j_1,\dots, j_t\}$. However, in the other cases, this condition is equivalent to $j_k -m \leq \sigma (j_k ) \leq m$ and 
$m+1 \leq \sigma(i_k) \leq m+i_k $ for all $k \in \gc 1,t \dc$. These observations imply that 
$\sigma \in \mathcal{S}_t$ if and only if there exist $m+1 \leq j_1 < \dots < j_t \leq 2m$ and $1 \leq i_1 < \dots < i_t \leq m$ with  
\begin{enumerate}
\item $j_k -m \leq \sigma (j_k ) \leq m \mbox{ for all } k \in \gc 1,t \dc$,
\item $ \sigma (l) \geq m+1 
\mbox{ for all } l \in \gc m+1,2m \dc \setminus \{j_1,\dots,j_t\} $,
\item $m+1 \leq \sigma(i_k) \leq m+i_k \mbox{ for all } k \in \gc 1,t \dc$
\item $ \mbox{ and } \sigma (l) \leq m 
\mbox{ for all } l \in \gc 1,m \dc \setminus \{i_1,\dots,i_t\} $.
\end{enumerate}
Hence we have 

\begin{eqnarray*}
\mid \mathcal{S}_t \mid & = & \left( \sum_{m+1 \leq j_1 < \dots < j_t \leq 2m} \left[ 2m-j_t+1 \right] \left[ 2m-j_{t-1} \right] \dots \left[ 2m-j_1+1-(t-1) \right] \right) \times (m-t)! \\
& & \times \left( \sum_{1 \leq i_1 < \dots < i_t \leq m} i_1 \left[ i_2 -1 \right] \dots \left[ i_t - (t-1) \right] \right)  \times (m-t)!
\end{eqnarray*}
Next, by means of the changes of variables $\alpha_k= 2m-j_{t+1-k }+1$, this leads to 
\begin{eqnarray*}
\mid \mathcal{S}_t \mid & = & \left( (m-t)! \sum_{1 \leq i_1 < \dots < i_t \leq m} i_1 \left[ i_2 -1 \right] \dots \left[ i_t - (t-1) \right] \right)^2. 
\end{eqnarray*}
Set $\displaystyle{ u(m,t):= \sum_{1 \leq i_1 < \dots < i_t \leq m} i_1 \left[ i_2 -1 \right] \dots \left[ i_t - (t-1) \right] }$ (with the convention that $u(m,0):=1$). In view of the previous equality, to prove that $\mid \mathcal{S}_t \mid =  \left[ (m-t)! S(m+1,m-t+1) \right]^2$, it just remains to prove that $u(m,t)= S(m+1,m-t+1)$ for all $m \geq 1$ and $t \in \gc 0, m \dc$. We proceed by induction on $m$.

If $m=1$, then obvious computations lead to $u(1,0)=1=S(2,2)$ and $u(1,1)=1=S(2,1)$, as desired. Assume now that $m > 1$. If $t=0$, then we have $u(m,0)=1=S(m+1,m+1)$, as desired. Let now $t \in \gc 1,m \dc$. Then we have 
$$u(m,t)= \sum_{1 \leq i_1 < \dots < i_{t-1} \leq m-1} i_1 \left[ i_2 -1 \right] \dots \left[ m - (t-1) \right]
+  \sum_{1 \leq i_1 < \dots < i_t \leq m-1} i_1 \left[ i_2 -1 \right] \dots \left[ i_t - (t-1) \right],$$ 
that is, $$u(m,t)= (m-t+1)u(m-1,t-1)+u(m-1,t).$$
So it follows from the induction hypothesis that $u(m,t)= (m-t+1)S(m,m+1-t)+S(m,m-t)$. Now it is well known that the Stirling numbers of second kind satisfy the following recurrence (see, for instance, \cite{stanleybook}):
$$S(m+1,m+1-t)=(m-t+1)S(m,m+1-t)+S(m,m-t),$$
so that $u(m,t)= S(m+1,m+1-t)$, as desired. \fin 
\\$ $



\section{Combinatorics of $\hc'_R$-$\spec(\sn)$.}
$ $

In this section, we recall the description of the set $\hc'_R$-$\spec(\sn)$ of those prime ideals of $\sn$ that are invariant under the natural action of the group $\hc'_R:=\{ (a_1,\dots,a_n,b_1,\dots,b_n) \in (\mathbb{C}^*)^{2n} \mid a_1 \dots a_n b_1 \dots b_n = 1 \}$ on this algebra. This description was obtained by Brown and Goodearl (see \cite{bg}) by using the partition of $\spec(\sn)$ constructed by Hodges and Levasseur (see \cite{hl1,hl2}) and will play a central role to describe the $\hc$-primes of $\mmpc$: we will construct in the next section an embedding from $\hc$-$\spec(\mmpc)$ into $\hc'_R$-$\spec(\sn)$ whose image will be determined in section \ref{sectionimagepsi} using this description of the set $\hc'_R$-$\spec(\sn)$.

\subsection{Torus actions on $\mnc$ and some related algebras.} 
\label{sectionnotations}
$ $

In this section, we fix some notations that will be used in the sequel of that paper. Recall that throughout this paper, we have fixed $n=m+p$. 
\\$ $
\\$\bullet$ We denote by $R$ the quantization of the ring of regular functions on 
$n \times n$ matrices with entries in $\mathbb{C}$; it is the $\mathbb{C}$-algebra generated by the 
$n \times n$ indeterminates $Y_{\ia}$, $1 \leq i,\alpha \leq n$, subject to the following 
relations.\\
If $\left( \begin{array}{cc} x & y \\ z & t \end{array} \right)$ is any $2 \times 2$ sub-matrix of 
$\mathcal{Y}:=\left( Y_{\ia} \right) _{(\ia) \in \gc 1,n \dc^2}$, then
 \begin{enumerate}
 \item $ yx=q^{-1} xy, \quad zx=q^{-1} xz, \quad zy=yz, \quad ty=q^{-1}
 yt, \quad tz=q^{-1} zt$.
 \item $tx=xt-(q-q^{-1})yz$.
\end{enumerate}
 It is well known that the ring $R$ is a Noetherian domain. Moreover, since 
$q$ is transcendental over $\mathbb{Q}$, it follows from \cite[Theorem 3.2]{gl4} that all prime ideals 
of $R$ are completely prime.
\\$\bullet$ It is well known that the group $\hc_R:=\left( \mathbb{C}^* \right)^{2n}$ acts 
on $R$ by $\mathbb{C}$-algebra automorphisms via:
$$(a_1,\dots,a_n,b_1,\dots,b_n).Y_{\ia} = a_i b_\alpha Y_{\ia} 
\quad \forall \: (\ia)\in \gc 1,n \dc^2.$$
An \textbf{$\hc_R $-eigenvector} $x$ of $R$ is a nonzero element $x \in R$ such that $h.x \in \mathbb{C}^*x$ 
for each $h \in \hc_R$. An ideal $I$ of $R$ is said to be \textbf{$\hc_R$-invariant} if $h.I =I$ for all $h\in \hc_R$. 
We denote by $\hc_R$-$\spec(R)$ the set of $\hc_R$-invariant prime ideals of $R$. Recall (see \cite[5.7. (i)]{gl1}) that $R$ has only finitely many $\hc_R$-invariant prime ideals.
\\$\bullet$ We now recall the notion of quantum minors. Let $K$ denote a $\mathbb{C}$-algebra. 
Let $u$ and $v$ be two positive integers and let $M=(x_{i,\alpha})_{(\ia) \in \gc 1,u \dc \times \gc 1,v \dc}$ be 
an $u \times v$ matrix with entries in $K$. The quantum minor of $M$ associated to the rows $i_1$, ..., $i_t$
 and the columns $\alpha_1$, ..., $\alpha_t$ ($t \in \gc 1, \min(u,v) \dc$, $1 \leq i_1 < \dots < i_t \leq u$, $1 \leq \alpha_1 < \dots < \alpha_t \leq v$) is the following element of $K$:
$$det_q(x_{\ia})_{\substack{i=i_1,\dots,i_t \\ \alpha=\alpha_1,\dots ,\alpha_t}}:=\sum_{\sigma \in S_t}(-q)^{l(\sigma)}x_{i_1,\alpha_{\sigma(1)}} \dots
x_{i_t,\alpha_{\sigma(t)}},$$
where $l(\sigma)$ denotes the length of the $t$-permuation $\sigma$.
\\$\bullet$ We call $\Delta$ the quantum determinant of $R$, that is, 
$$\Delta :=det_q(\yia)_{(\ia) \in \gc 1,n \dc^2}.$$
Recall that $\Delta$ belongs to the center of $R$. 
\\$\bullet$ We denote by $\gn$ the quantization of the ring of regular functions on 
$GL_n(\mathbb{C})$, that is, $\gn:=R[\Delta^{-1}]$. Recall that $\gn$ is a Noetherian domain. Moreover, since $q$ is transcendental over $\mathbb{Q}$, it follows from \cite[Corollary II.6.10]{bg2} that all prime ideals 
of $\gn$ are completely prime. Note also that, since $\Delta $ is an $\hc_R$-eigenvector of $R$, the action of $\hc_R$ 
on $R$ induces an action of $\hc_R$ on $\gn$.
\\$\bullet$ Let $\sn$ be the quantization of the ring of regular functions on $SL_n(\mathbb{C})$ and by $X_{\ia}$ ($(\ia) \in \gc 1,n \dc^2$) the canonical generators of $\sn$, that is, we set
$$\sn := \frac{R}{\left\langle \Delta -1 \right\rangle } \mbox{ and } X_{\ia} := \yia +\left\langle \Delta -1 \right\rangle \mbox{ for all }(\ia) \in \gc 1,n \dc^2. $$ 
Recall that $\sn$ is a Noetherian domain (see \cite{ls}). Moreover, since $q$ is transcendental over $\mathbb{Q}$, it follows from \cite[Corollary II.6.10]{bg2} that all prime ideals 
of $\sn$ are completely prime.
\\$\bullet$ Note that, because of the relation $\Delta=1$, the action of the torus $\hc_R$ on $R$ does not induce a corresponding action of $\hc_R$ on $\sn$. To deal with this problem, we consider the stabilizer $\hc'_R$ of the quantum determinant, that is, we set  
$$\hc'_R:=\{ (a_1,\dots,a_n,b_1,\dots,b_n) \in \hc_R \mid a_1 \dots a_n b_1 \dots b_n = 1 \}.$$
This group $\hc'_R$ stabilizes $\Delta$ and so acts on $\sn$ by $\mathbb{C}$-algebra automorphisms via:
$$(a_1,\dots,a_n,b_1,\dots,b_n).X_{\ia} = a_i b_\alpha X_{\ia} 
\quad \forall \: (\ia)\in \gc 1,n \dc^2.$$
An \textbf{$\hc'_R $-eigenvector} $x$ of $\sn$ is a nonzero element $x \in \sn$ such that $h.x \in \mathbb{C}^*x$ 
for each $h \in \hc'_R$. An ideal $I$ of $\sn$ is said to be \textbf{$\hc'_R$-invariant} if $h.I =I$ for all $h\in \hc'_R$. 
We denote by $\hc'_R$-$\spec(\sn)$ the set of $\hc'_R$-invariant prime ideals of $\sn$. Since $q$ is transcendental, it follows from \cite[Theorem II.5.17]{bg2} that $\hc'_R$-$\spec(\sn)$ is a finite set.
\\$ $

\subsection{Description of $\hc'_R$-$\spec(\sn)$.}
$ $

In this section we recall the description of the set $\hc'_R$-$\spec(\sn)$ that was obtained by Brown and Goodearl (see \cite{bg}).

\begin{nota}
\label{notationsideauxhodgeslevasseur}
$ $
\begin{itemize}
\item For all $j\in \gc 1,n-1 \dc$ and $y\in S_n$, we set $c_{j,y}^+:=det_q(\xia)_{(\ia) \in y(\gc 1, j \dc ) \times \gc 1,j \dc}$ and $c_{j,y}^-:=det_q(\xia)_{(\ia) \in y(\gc j+1,n \dc ) \times \gc j+1,n \dc}$.
\item Let $w=(w_+,w_-) \in S_n \times S_n$. We denote by $I_{w^+}^+$ the ideal of $\sn$ generated by those $c_{j,y}^{+}$ with $j\in \gc 1,n-1 \dc $, $y \in S_n$ and $y \nleq_j w_{+}$, that is, 
$$I_{w_+}^+:= \langle c_{j,y}^{+} \mid  j\in \gc 1,n-1 \dc \mbox{, }y \in S_n \mbox{ and }y \nleq_j w_{+} \rangle.$$
We denote by $I_{w^-}^-$ the ideal of $\sn$ generated by those $c_{j,y}^{-}$ with $j\in \gc 1,n-1 \dc $, $y \in S_n$ and $y \nleq_j w_{-}$, that is, 
$$I_{w_-}^-:= \langle c_{j,y}^{-} \mid  j\in \gc 1,n-1 \dc \mbox{, }y \in S_n \mbox{ and }y \nleq_j w_{-} \rangle,$$
and we set $I_w:=I_{w_+}^+ + I_{w_-}^-$.
\end{itemize}
\end{nota}

Observe that, since $I_w$ is generated by quantum minors, $I_w$ is an $\hc'_R$-invariant ideal of $\sn$. 
Further, since $q$ is transcendental, it follows from \cite[Theorem 3]{jos2} that $I_w$ is (completely) prime. Hence each $I_w$ is an $\hc'_R$-invariant prime ideal of $\sn$. In fact, there is no other $\hc'_R$-invariant prime ideal in $\sn$ since Brown and Goodearl have shown (see \cite{bg}; see also \cite[Proposition 3.4.6]{lau}) that 

\begin{prop}
\label{propHprimesln}
$\hc'_R$-$\spec(\sn) = \{ I_w \mid w \in S_n \times S_n \}$.
\end{prop}

In the sequel, $S_n \times S_n$ will always be endowed with the product order, that is, $w=(w_+,w_-) \leq \theta=(\theta_+,\theta_-)$ if and only if $w_+ \leq \theta_+$ and $w_- \leq \theta_-$ (where $\leq$ still denotes the reverse Bruhat order on $S_n$).

\begin{prop}
\label{bijsn}
The map $\chi :   S_n \times S_n  \rightarrow  \hc'_R \mbox{-}\spec(\sn)$ 
defined by $\chi(w)=I_w$ is an order-reversing bijection; its inverse is also an order-reversing bijection.
\end{prop}
\preuve In view of Proposition \ref{propHprimesln}, it just remains to show that $\chi$ and $\chi^{-1}$ are decreasing. Concerning $\chi$, it easily follows from the characterization of the (reverse) Bruhat ordering given in Proposition \ref{propbruhat} and the definition of $I_w$. So it just remains to deal with $\chi^{-1}$. Let $J \subseteq K$ be two $\hc'_R$-invariant prime ideals of $\sn$. There exist $w_J$ and $w_K$ in $S_n \times S_n$ such that $J=I_{w_J}$ and $K=I_{w_K}$. In order to prove that $\chi^{-1}$ is decreasing, we need to show that 
$w_J = (w_J^+,w_J^-) \geq w_K = (w_K^+,w_K^-)$, that is, $w_J^+ \geq w_K^+$ and $w_J^- \geq w_K^-$. Assume that this is not the case. Then we have for instance $w_J^+ \ngeq w_K^+$. It follows from Proposition \ref{propbruhat} that there exists $j \in \gc 1,n-1 \dc$ such that $w_J^+ \ngeq_j w_K^+$ and so we have $c_{j,w_K^+}^+ \in I_{w_J}=J \subseteq K$. Hence the quantum minor $c_{j,w_K^+}^+$ must belong to $K=I_{w_K}$. But, on the other hand, Hodges and Levasseur have shown (see \cite[Theorem 2.2.3]{hl1}) that, for all $i \in \gc 1,n-1 \dc$, $c_{i,w_K^+}^+ \not \in I_{w_K} $. This is a contradiction and so we have $w_J^+ \geq w_K^+$ and $w_J^- \geq w_K^-$, as required. \fin


\section{An order-preserving embedding from $\hc$-$\spec(\mmpc)$ into $\hc'_R$-$\spec(\sn)$.}
$ $

Recall that $A= \mmpc$. Following the route sketched in the introduction, we construct in this section an order-preserving embedding $\psi : \hc \mbox{-} \spec(A) \longrightarrow \hc'_R \mbox{-} \spec(\sn)$. In order to do this, we will use the "$(m,p)$ deleting-derivations algorithm" introduced in \cite{lau}. This algorithm consists of certain changes of variables in the field of fractions of $R=\mnc$. At the step $(m,m)$ of this algorithm, the subalgebra $R^{(m,m)}$ of $\fract(R)$ generated by the new indeterminates is isomorphic to an iterated Ore extension of $A$ that does not involve $\sigma$-derivations. This fact allows the construction of an embedding from $\hc$-$\spec(A)$ into a subset of $\spec(R^{(m,m)})$, the image of $J \in \hc$-$\spec(A)$ just being the ideal of $R^{(m,m)}$ generated by $J$. Next, using some suitable localizations and contractions related to the involved changes of variables, we are able to extend this embedding to an embedding from $\hc$-$\spec(A)$ into $\{L \in \hc_R \mbox{-} \spec(R) \mid \Delta \not \in P \}$. The end of the construction of the embedding $\psi$ uses some classical results from localization theory and the isomorphim $\sn[z,z^{-1}] \simeq \gn$ that was constructed by Levasseur and Stafford in \cite{ls}. Note that the construction of this embedding $\psi$ is implicit in \cite{lau} and so, in this section, we will essentially recall results of \cite{lau}.

Throughout this section, we keep the notations and conventions of section \ref{sectionnotations}.
\\$ $ 

\subsection{The $(m,p)$ deleting-derivations algorithm in $R$.}
$ $


Define the relation "$\leq_m $" by
$$(\ia) \leq_m (j,\beta) \Longleftrightarrow 
\left\{ \begin{array}{ll}
\left[ ( i < j ) \mbox{ or } ( i = j \mbox{ and } \alpha \leq \beta ) \right] &  \mbox{ if } j >m \\
 \left[ ( i < m) \mbox{ and } \left( ( \alpha > \beta ) \mbox{ or } ( \alpha = \beta \mbox{ and } i \geq j ) \right) \right] &  \mbox{ if } j \leq m  
\end{array}
 \right.$$
This defines a total ordering on $\mathbb{N}^2$ that we call the \textit{$(m,p)$-ordering on $\mathbb{N}^2$}.\\

\begin{nota}
\label{notareferee}
$ $
\begin{itemize}
\item We set $E=\left( \gc 1,n \dc^2 \cup \{(n,n+1)\} \right) \setminus \{(m,n)\}$.
\item Let $(j,\beta) \in E \cup \{(m,n)\}$. If $(j,\beta) \neq (n,n+1)$, then we denote by $(j,\beta)^{+}$ the smallest element (relative to $\leq_m$) of the set $\{(\ia) \in E \mid (\ia) >_m (j,\beta) \}$.
\end{itemize}
 \end{nota}

Note that, for all $(\ia ) \in E$, we have $(m-1,n) \leq_m (\ia) \leq_m  (n,n+1)$.\\$ $

Recall (see \cite[Proposition 3.2.6]{lau}) that $R$ can be presented as an iterated Ore extension over $\mathbb{C}$, with the generators adjoined in the $(m,p)$-ordering. Moreover \cite[Proposition 3.2.6]{lau} shows that the theory of the deleting-derivations (see \cite{c1}) can be applied to this iterated Ore extension. The corresponding algorithm is called the $(m,p)$ deleting-derivations algorithm. It is known (see \cite[paragraph 3]{lau}) that this algorithm consists of the construction, for each $r \in E$, of a family $(Y_{\ia}^{(r)})_{(\ia) \in \gc 1,n \dc^2 }$ of elements of $F=\fract(R)$, defined as follows.
\begin{enumerate}
\item If $r=(n,n+1)$, then $Y_{\ia}^{(r)}=Y_{\ia}$ for all $(\ia) \in \gc 1,n \dc^2 $.
\item Assume that $r=(j,\beta) <_m (n,n+1)$ and that the $Y_{\ia}^{(r^{+})}$ ($(\ia) \in \gc 1,n \dc^2 $) are already known. If $(\ia) \in \gc 1,n \dc^2 $, then 
$$Y_{\ia}^{(r)}=
\left\{ \begin{array}{ll}
Y_{\ia}^{(r^{+})}-Y_{i,\beta}^{(r^{+})}\left(Y_{j,\beta}^{(r^{+})}\right)^{-1}
 Y_{j,\alpha}^{(r^{+})} 
& \mbox{ if }  j < i \leq m \mbox{ and }  \beta < \alpha \\
 
 Y_{\ia}^{(r^{+})}-Y_{i,\beta}^{(r^{+})}\left(Y_{j,\beta}^{(r^{+})}\right)^{-1}
 Y_{j,\alpha}^{(r^{+})}  
& \mbox{ if } i < j \mbox{, } \alpha < \beta \mbox{ and }j >m \\

Y_{\ia}^{(r^{+})} & \mbox{ otherwise.}

\end{array} \right.$$
\end{enumerate} 
$ $

Observe that the notations of this paper are not exactly the same as those of \cite{lau}. For convenience, we have dropped throughout this paper almost every subscript $m$. \\$ $

As in \cite{lau}, if $r\in E$, we denote by $R^{(r)}$ the subalgebra of $F=\fract(R)$ generated by the $Y_{\ia}^{(r)}$ ($(\ia) \in \gc 1,n \dc^2 $), that is,
$$R^{(r)}:=\mathbb{C} \langle Y_{\ia}^{(r)} \mbox{ $\mid$ }(\ia) \in \gc 1,n \dc^2  \rangle .$$
Moreover, we set $\ov{R}:=R^{(m-1,n)}$ and $V_{\ia}:=Y_{\ia}^{(m-1,n)}$ for all $(\ia) \in \gc 1,n \dc^2$.
\\$ $

The following observation is an easy consequence of the above formulae (that is, those formulae that express the $Y_{\ia}^{(j,\beta)}$ in terms of the $Y_{\ia}^{(j,\beta)^+}$). \\

\begin{obs}
Let $(\ia) \in \gc 1,n \dc^2$ and $(j,\beta) \in E$ with $(j,\beta) \leq_m (i,\alpha)^+$.
\\We have $V_{\ia}= Y_{\ia}^{(j,\beta)}$, so that $V_{\ia}$ belongs to $R^{(j,\beta)}$. 
\end{obs}
$ $

Let $(j,\beta) \in E$. It follows from the theory of deleting-derivations that $R^{(j,\beta)}$ is a Noetherian domain (see \cite[Th\'eor\`eme 3.2.1]{c1}) and that $\fract(R^{(j,\beta)})=\fract(R)$ 
(see \cite[Th\'eor\`eme 3.3.1]{c1}). Moreover (see \cite[section 3.2]{lau}) all prime ideals of $R^{(j,\beta)}$ are completely prime and the torus $\hc_R$ acts by automorphisms on $R^{(j,\beta)}$ via:
$$(a_1,\dots,a_n,b_1,\dots,b_n).\yia^{(j,\beta)} =a_ib_{\alpha} \yia^{(j,\beta)} $$
for all $(\ia) \in \gc 1,n \dc ^2$ and $(a_1,\dots,a_n,b_1,\dots,b_n) \in \hc_R $.\\

Finally, recall (see \cite[Lemme 3.2.11]{lau}) that the quantum determinant $\Delta$ of $R$ belongs to all these algebras, that is, 
$$\Delta \in R^{(j,\beta)}, \: \forall (j,\beta) \in E,$$
and that, in $R^{(m,m)}$, $\Delta$ is just the product of the diagonal indeterminates $V_{i,i}=Y_{i,i}^{(m,m)}$, that is,
$$\Delta =Y_{1,1}^{(m,m)} Y_{2,2}^{(m,m)}\dots Y_{n,n}^{(m,m)}=V_{1,1} V_{2,2} \dots V_{n,n}.$$
$ $

\subsection{An order-preserving embedding $\iota : \hc \mbox{-}\spec(A) \rightarrow \hc_R \mbox{-}\spec^{**}(R^{(m,m)})$.}
$ $

One of the reasons that makes the $(m,p)$ deleting-derivations algorithm interesting is that it follows from the theory of the deleting-derivations that there is a natural embedding from $A$ into the algebra $R^{(m,m)}$ obtained from $R$ at the step $(m,m)$ of the algorithm. 
More precisely, the subalgebra of $R^{(m,m)}$ generated by the $\yia^{(m,m)}$ with $(\ia) \in \gc 1,m \dc \times \gc m+1,n \dc$ can be identified to $A$, and the $\yia^{(m,m)}$ with $(\ia) \in \gc 1,m \dc \times \gc m+1,n \dc$ to the canonical generators of $A$ (see \cite[3.3]{lau}). Thus we denote by $A$ this subalgebra of $R^{(m,m)}$. Moreover $R^{(m,m)}$ can be written as an iterated Ore extension of $A$ that does not involve 
$\sigma$-derivations. This fact has allowed us to prove the following statement (see \cite[Proposition 3.3.3]{lau}). \\

\begin{prop}
\label{propextensionMR}
 Let $J$ be an $\hc$-invariant prime ideal in $A$. 
Set 
$$I:=\sum_{a_{\ia} \in \mathbb{N} } J V_{m,m}^{a_{m,m}} V_{m-1,m}^{a_{m-1,m}}\dots V_{1,1}^{a_{1,1}} V_{m+1,1}^{a_{m+1,1}} \dots  V_{n,n}^{a_{n,n}}.$$
$I$ is an $\hc_R$-invariant prime ideal of $R^{(m,m)}$ such that 
\begin{enumerate}
\item $I \cap A = J$ and
\item $V_{\ia}  \notin I$ for all $(m,m) \leq_m (\ia) \leq_m (n,n)$.
\end{enumerate}
\end{prop} 
 $ $

As in \cite[Conventions 3.3.4]{lau}, we need to introduce the following notations.\\

\begin{nota}
\label{convinjcano}
Let $(j,\beta) \in E$ with $(j,\beta) \neq (n,n+1)$.
\\$\bullet$ We set 
$$\spec^* \left( R^{(j,\beta)} \right) := \{ P \in \spec \left( R^{(j,\beta)} \right) \mid 
(\forall u \in \gc 1,n \dc^2 ) 
\mbox{ }((j,\beta) \leq_m u \mbox{ } \Rightarrow \mbox{ } V_u \notin P) \}$$ 
and 
$$\hc_R \mbox{-}\spec^* \left( R^{(j,\beta)} \right) := \hc_R \mbox{-}\spec \left( R^{(j,\beta)} \right) \cap
\spec^* \left( R^{(j,\beta)} \right).$$
\\$\bullet$ We set $$\spec^{**} \left( R^{(j,\beta)} \right) := \{ P \in \spec^* \left( R^{(j,\beta)} \right) \mid 
 \Delta \notin P \}$$ and 
$$\hc_R \mbox{-}\spec^{**} \left( R^{(j,\beta)} \right) := \hc_R \mbox{-}\spec \left( R^{(j,\beta)} \right) \cap
\spec^{**} \left( R^{(j,\beta)} \right).$$
\end{nota}
$ $

Let $J$ be an $\hc$-invariant prime ideal in $A$. 
Set 
$$I:=\sum_{a_{\ia} \in \mathbb{N} } J V_{m,m}^{a_{m,m}} V_{m-1,m}^{a_{m-1,m}}\dots V_{1,1}^{a_{1,1}} V_{m+1,1}^{a_{m+1,1}} \dots  V_{n,n}^{a_{n,n}}.$$
Since $I$ is a completely prime ideal of $R^{(m,m)}$ which does not contain the 
elements $Y_{k,k}^{(m,m)}=V_{k,k}$ $(k \in \gc 1,n \dc$) (see Proposition \ref{propextensionMR}), $I$ does not contain the product $V_{1,1} V_{2,2} \dots V_{n,n}$. Since this product is equal to the quantum determinant $\Delta$ (see \cite[Lemme 3.2.11]{lau}), we obtain that $\Delta$ does not belong to $I$. On the other hand, it follows from Proposition \ref{propextensionMR} that $I \in 
\hc_R \mbox{-}\spec^* \left( R^{(j,\beta)} \right)$. Hence we conclude that 
$I \in \hc_R \mbox{-}\spec^{**} \left( R^{(j,\beta)} \right)$. This allows us to define a map 
$\iota : \hc \mbox{-}\spec(A) \rightarrow \hc_R \mbox{-}\spec^{**} (R^{(m,m)})$ by
$$\iota (J) = \sum_{a_{\ia} \in \mathbb{N} } J V_{m,m}^{a_{m,m}}V_{m-1,m}^{a_{m-1,m}} \dots V_{1,1}^{a_{1,1}} V_{m+1,1}^{a_{m+1,1}} \dots  V_{n,n}^{a_{n,n}}. $$
An immediate consequence of Proposition \ref{propextensionMR} is the following result.\\

\begin{prop}
\label{iota}
The map $\iota : \hc \mbox{-}\spec(A) \rightarrow \hc_R \mbox{-}\spec^{**} (R^{(m,m)})$ 
 defined by
$$\iota (J) = \sum_{a_{\ia} \in \mathbb{N} } J V_{m,m}^{a_{m,m}}V_{m-1,m}^{a_{m-1,m}} \dots V_{1,1}^{a_{1,1}} V_{m+1,1}^{a_{m+1,1}} \dots  V_{n,n}^{a_{n,n}} $$
is an ordered embedding.
\end{prop}
 $ $

\subsection{An ordered embedding $\tau : \hc_R \mbox{-} \spec^{**} ( R^{(m,m)} ) \rightarrow \{ P \in \hc_R \mbox{-} \spec ( R ) \mid \Delta \not \in P \}$.}
\label{sectionLjbeta}
$ $

In \cite[3.3]{lau}, we have associated to an ideal $I \in \hc_R \mbox{-} \spec^{**} ( R^{(m,m)} )$ a unique ideal $L \in \{ P \in \hc_R \mbox{-} \spec ( R ) \mid \Delta \not \in P \}$ by using some suitable localizations and contractions related to the $(m,p)$ deleting-derivations algorithm. In this section, we show that this construction leads to an ordered embedding $\tau : \hc_R \mbox{-} \spec^{**} ( R^{(m,m)} ) \rightarrow \{ P \in \hc_R \mbox{-} \spec ( R ) \mid \Delta \not \in P \}$ .

\begin{nota}
 We denote by $\Sigma_{j,\beta}$ the multiplicative system of $R^{(j,\beta)}$ (resp. $R^{(j,\beta)^+}$) generated by $V_{j,\beta}=Y_{j,\beta}^{(j,\beta)}=Y_{j,\beta}^{(j,\beta)^+}$. 
\end{nota}
$ $

Recall (see \cite[Th\'eor\`eme 3.2.1]{c1}) that, for all $(j,\beta) \in E$ with $(j,\beta) \neq (n,n+1)$, $\Sigma_{j,\beta}$ is a denominator set in both $R^{(j,\beta)}$ and $R^{(j,\beta)^+}$, and that 
$$R^{(j,\beta)} \Sigma_{j,\beta}^{-1}  = R^{(j,\beta)^+} \Sigma_{j,\beta}^{-1} .$$

Let $(j,\beta) \in E$ with $ (m,m) \leq_m (j,\beta) \leq_m (n,n)$ and $P \in \hc_R \mbox{-} \spec^{**} \left( R^{(j,\beta)} \right)$. Then $P \cap  \Sigma_{j,\beta} = \emptyset$, so that the ideal $P\Sigma_{j,\beta}^{-1}$ is a completely prime ideal of $R^{(j,\beta)} \Sigma_{j,\beta}^{-1}  = R^{(j,\beta)^+} \Sigma_{j,\beta}^{-1}$. Hence $P\Sigma_{j,\beta}^{-1} \cap R^{(j,\beta)^+}$ is a completely prime ideal of $R^{(j,\beta)^+}$. Moreover, since $P\in \hc_R \mbox{-} \spec^{**} \left( R^{(j,\beta)} \right)$, we deduce from 
\cite[Lemme 3.3.5 and Lemme 3.3.8]{lau} that $P\Sigma_{j,\beta}^{-1} \cap R^{(j,\beta)^+} \in  \hc_R \mbox{-} \spec^{**} \left( R^{(j,\beta)^+} \right)$. This allows us to define, for all $(j,\beta) \in E$ with $ (m,m) \leq_m (j,\beta) \leq_m (n,n)$, a map $\tau_{j,\beta} : \hc_R \mbox{-} \spec^{**} \left( R^{(j,\beta)} \right) \rightarrow \hc_R \mbox{-} \spec^{**} \left( R^{(j,\beta)^+} \right)$ by $\tau_{j,\beta} (P) = P \Sigma_{j,\beta}^{-1} \cap R^{(j,\beta)^+}$. Moreover, it is easy to show that

\begin{prop}
\label{taujbeta}
For all $(j,\beta) \in E$ with $ (m,m) \leq_m (j,\beta) \leq_m (n,n)$, the map $\tau_{j,\beta} : \hc_R \mbox{-} \spec^{**} \left( R^{(j,\beta)} \right) \rightarrow \hc_R \mbox{-} \spec^{**} \left( R^{(j,\beta)^+} \right)$ defined by $\tau_{j,\beta} (P) = P \Sigma_{j,\beta}^{-1} \cap R^{(j,\beta)^+}$ is an order-preserving embedding.
\end{prop}
$ $

Let $(j,\beta) \in E$ with $ (m,m) \leq_m (j,\beta) \leq_m (n,n)$ and $P \in \hc_R \mbox{-} \spec^{**} \left( R^{(j,\beta)} \right)$. Note that we can express $P$ in terms of $\tau_{j,\beta} (P)$ as follows.\\$ $

\begin{obs}
\label{obsutile} 
$V_{j,\beta} \not \in \tau_{j,\beta} (P)$ and $P = \tau_{j,\beta} (P)\Sigma_{j,\beta}^{-1} \cap R^{(j,\beta)}$. 
\\In other words, $P= \varphi_{j,\beta} (\tau_{j,\beta}(P))$, where 
$\varphi_{(j,\beta)}:\{Q \in \spec(R^{(j,\beta)^+}) \mid V_{j,\beta} \not \in Q \} \rightarrow
\spec(R^{(j,\beta)})$ denotes the map defined by $\varphi_{j,\beta} (Q) := Q \Sigma_{j,\beta}^{-1} \cap R^{(j,\beta)}$.
\end{obs}
 
Composing the ordered embeddings $\tau_{\ia}$ leads to the following statement.

\begin{prop}
\label{tau}
The composition map 
$\tau := \tau_{n,n} \circ \dots \circ  \tau_{m+1,1} \circ \tau_{1,1} \circ \dots \circ \tau_{m,m} $ is an ordered embedding from $ \hc_R \mbox{-} \spec^{**} \left( R^{(m,m)} \right)$ into $\{ P \in \hc_R \mbox{-} \spec \left( R \right) \mid \Delta \notin P\}$.
\end{prop}
$ $




\subsection{An ordered bijection between $\{ P \in \hc_R \mbox{-} \spec \left( R \right) \mid \Delta \not \in P\}$ and $\hc'_R$-$\spec(\sn)$.}
\label{sectionSn}
$ $

In this section, we deduce from classical results of localization theory and the isomorphism  
$\sn[z,z^{-1}] \simeq \gn$ exhibited in \cite{ls} the construction of an order-preserving bijection 
between $\{ P \in \hc_R \mbox{-} \spec \left( R \right) \mid \Delta \notin P\}$ and $\hc'_R$-$\spec(\sn)$.
\\$ $

First it is well known (see \cite[Exercise II.1.J]{bg2}) that the extension 
$\epsilon :  \{ P \in \hc_R \mbox{-} \spec \left( R \right) \mid \Delta \notin P\} \rightarrow \hc_R \mbox{-} \spec (\gn) = \hc_R \mbox{-} \spec (R[\Delta^{-1}])$ defined by $ \epsilon (P)=P_{\Delta}=P[\Delta^{-1}]$ is an ordered bijection whose inverse is just the contraction.

Next recall that Levasseur and Stafford (see \cite{ls}) have constructed an isomorphism between $\sn[z,z^{-1}]$ and $\gn$   from which Brown and Goodearl have deduced the existence of an ordered bijection 
$\varphi : \hc'_R \mbox{-}\spec(\sn) \rightarrow \hc_R \mbox{-} \spec(\gn)$ . Let us now make precise the construction of this bijection.

\begin{prop}[see \cite{ls}]
There exists an algebra isomorphism $\theta : \sn[z,z^{-1}] \rightarrow \gn$ with
$$\left\{ \begin{array}{llll}
\theta (X_{\ia}) & = & Y_{\ia} & \mbox{ if } i>1 \\
\theta (X_{1,\alpha}) & = & Y_{1,\alpha} \Delta^{-1} & \\
\theta (z) & = & \Delta. & 
\end{array} \right.$$
\end{prop}

\begin{nota}
Let $P  \in \spec(\sn)$. We set $\widehat{P}=\displaystyle{\bigoplus_{i \in \mathbb{Z}}P z^i}$.
\end{nota}

Note that, for all $P  \in \spec(\sn)$, $\widehat{P}$ is a (completely) prime ideal of $\sn[z,z^{-1}]$ and that $\widehat{P} \cap \sn=P$.

In \cite[Lemma II.5.16]{bg2}, Brown and Goodearl have established the following result.\\

\begin{lem}
\label{rappelbijHpremier}
The map $\varphi :  \hc'_R \mbox{-}\spec(\sn) \rightarrow   \hc_R \mbox{-} \spec(\gn) $ defined by 
$\varphi(P)= \theta(\widehat{P})$ is an ordered bijection. Moreover its inverse $\varphi^{-1} :   \hc_R \mbox{-} \spec(\gn) \rightarrow \hc'_R \mbox{-}\spec(\sn)$ is also an ordered bijection.
\end{lem}

Now composing the map $\varphi^{-1}$ and $\epsilon$ we get a ordered bijection from $\{ P \in \hc_R \mbox{-} \spec \left( R \right) \mid \Delta \notin P\}$ onto $\hc'_R \mbox{-}\spec(\sn)$, whose inverse is also ordered.\\

\subsection{An ordered embedding $\psi : \hc \mbox{-} \spec( A ) \rightarrow \hc'_R \mbox{-} \spec (\sn)$.}
\label{sectionpsi}
$ $

By composing the maps $\iota$, $\tau $, $\epsilon$ and $\varphi^{-1}$ of the previous sections, we easily obtain the following statement.\\

\begin{prop} 
The map $\psi :=\varphi^{-1} \circ \epsilon \circ \tau \circ \iota : \hc \mbox{-} \spec (A) \rightarrow \hc'_R \mbox{-} \spec (\sn )$ is an order-preserving embedding.
\end{prop}
$ $

Naturally, $\psi$ induces an ordered bijection $\widetilde{\psi}$ from  $\hc \mbox{-} \spec (A)$ onto $\psi \left( \hc \mbox{-} \spec (A) \right)$. Observe that, because of the construction of $\psi$ (see Proposition \ref{propextensionMR} and Observation \ref{obsutile}), it is easy to check that the inverse of $\widetilde{\psi}$ is also an ordered bijection. Moreover, in view of Proposition \ref{propextensionMR}, Observation \ref{obsutile} and the results of the previous paragraph, it is easy to check that, for all $K \in  \psi \left( \hc \mbox{-} \spec (A) \right)$, we have
$$\widetilde{\psi}^{-1}(K)= \psi_{(m,m)}\left( \varphi(K) \cap R \right) \cap A,$$
where $\psi_{(m,m)}:=\varphi_{(m,m)}\circ \varphi_{(m-1,m)} \circ \dots \circ \varphi_{(n,n)}$ and $\varphi_{(j,\beta)}:\{Q \in \spec(R^{(j,\beta)^+}) \mid V_{j,\beta} \not \in Q \} \longrightarrow
\spec(R^{(j,\beta)})$ denotes the map defined by $\varphi_{j,\beta} (Q) := Q \Sigma_{j,\beta}^{-1} \cap R^{(j,\beta)}$.


\section{Effect of the $(m,p)$ deleting-derivations algorithm on quantum minors.}
$ $

In the previous section \ref{sectionpsi}, we have constructed an order-preserving embedding $\psi : \hc \mbox{-} \spec (A) \rightarrow \hc'_R \mbox{-} \spec (\sn )$. Naturally we want to determine its image. Recall that Brown and Goodearl have determined the set $\hc'_R \mbox{-} \spec (\sn )$. More precisely, recall that they have shown (see Proposition \ref{propHprimesln}) that $\hc'_R \mbox{-} \spec (\sn ) = \{ I_w \mid w \in S_n \times S_n \}$, where each $I_w$ is generated by some quantum minors. So, in order to determine the image of $\psi$, that is, which $I_w$ belongs to the image of $\psi$, we need first to obtain some technical criteria for a quantum minor to belong to $\psi(J)$ ($J \in \hc $-$\spec(A)$). In this section, we establish such criteria. In the following section, we will use these new tools 
in order to prove that the image of $\psi$ is exactly the set of those $I_{w_0,w_0\sigma}$ such that $\sigma \in  \mathcal{S}$ (where $w_0$ still denotes the longest element of $S_n$ and $\mathcal{S}$ the sub-poset of the (reverse) Bruhat order of $S_n$ that we have studied in section \ref{sectionbruhatS}).

\subsection{Quantum minors and $(m,p)$ deleting-derivations algorithm.} 
$ $

We keep the conventions and notations of the previous sections. Moreover, if $(j,\beta) \in \gc 1,n \dc ^2$, then
 $R_{j,\beta}^{(j,\beta)^+}$ denotes the subalgebra of $\fract(R)$ generated by the $Y_{i,\alpha}^{(j,\beta)^+}$ with 
$(m,n) \leq_m (i,\alpha) <_m (j,\beta)$, that is, 
$$R_{j,\beta}^{(j,\beta)^+} : = \mathbb{C} \langle Y_{i,\alpha}^{(j,\beta)^+} \mid (m,n) \leq_m (i,\alpha) <_m (j,\beta) \rangle ,$$
 and $R_{j,\beta}^{(j,\beta)}$ denotes the subalgebra of $\fract(R)$ generated by the $Y_{i,\alpha}^{(j,\beta)}$ with $(m,n) \leq_m (i,\alpha) <_m (j,\beta)$, that is, 
$$R_{j,\beta}^{(j,\beta)} : = \mathbb{C} \langle Y_{i,\alpha}^{(j,\beta)} \mid (m,n) \leq_m (i,\alpha) <_m (j,\beta) \rangle .$$
It follows from \cite[Th\'eor\`eme 3.2.1]{c1} that, for all $(j,\beta) \in  \gc 1,n \dc ^2$, there exists a (unique) $\mathbb{C}$-algebra homomorphism $\theta_{j,\beta} : R_{j,\beta}^{(j,\beta)^+} \rightarrow R_{j,\beta}^{(j,\beta)}$ which sends  $Y_{i,\alpha}^{(j,\beta)^+}$ to $Y_{\ia}^{(j,\beta)}$ ($(m,n) \leq_m (i,\alpha) <_m (j,\beta)$).

Let now $J$ be a fixed $\hc$-invariant prime ideal of $A$ and set $L:= \tau \circ \iota (J)$; this is an $\hc_R$-invariant prime ideal of $R$ that does not contain the quantum determinant $\Delta$. We set $L^{(m,m)}:= \iota (J)$;  further, for all $(j,\beta) \in E$ with $(m-1,m) \leq_m (j,\beta) \leq_m (n,n+1)$, we set $L^{(j,\beta)}:= \tau_{(j,\beta)^-} \circ \dots \circ \tau_{m,m} \circ \iota (J)$ where $(j,\beta)^-$ denotes the greatest element (relative to $\leq_m$) of the set $\{(\ia) \in E \mid (\ia) < (j,\beta)\}$; $L^{(j,\beta) }$ is an $\hc_R$-invariant prime ideal of $R^{(j,\beta)}$.

Recall (see Observation \ref{obsutile}) that, for all $(j,\beta) \in E$ with $(m,m) \leq_m (j,\beta) \leq_m (n,n)$, we have $L^{(j,\beta)^+}= L^{(j,\beta)} \Sigma_{j,\beta}^{-1} \cap R^{(j,\beta)^+}$ and $L^{(j,\beta)}= L^{(j,\beta)^+} \Sigma_{j,\beta}^{-1} \cap R^{(j,\beta)}$.

\begin{nota}
$ $
\\$\bullet$ We set $\displaystyle{B:=\frac{R}{L}}$; $B$ is a Noetherian domain. We denote by $G$ its skew-field of fractions.
\\$\bullet$ For all $(j,\beta) \in E$ with $(m,m) \leq_m (j,\beta) \leq_m (n,n+1)$, 
we set $\displaystyle{B^{(j,\beta)}:=\frac{R^{(j,\beta)}}{L ^{(j,\beta)}}}$ 
and $y_{i,\alpha}^{(j,\beta)}:=Y_{i,\alpha}^{(j,\beta)}+L ^{(j,\beta)}$ for each $(i,\alpha) \in \gc 1,n \dc^2$. 
\end{nota}
$ $

Note that, for all $(j,\beta) \in E$, the ring $B^{(j,\beta)}$ is a Noetherian domain whose skew-field of fractions is $G$.

\begin{nota}
For all $(j,\beta),(k,\gamma) \in E$ with $(m,m) \leq_m (j,\beta) \leq_m (n,n+1)$, we denote by $B_{k,\gamma}^{(j,\beta)}$ the subalgebra of $B^{(j,\beta)}$ defined by
$$B_{k,\gamma}^{(j,\beta)}:=\mathbb{C} \langle y_{i,\alpha}^{(j,\beta)} \mbox{ $\mid$ }(m,n)  \leq_m (i,\alpha)  <_m  (k,\gamma) \rangle.$$
\end{nota}
$ $

Let $(j,\beta) \in E$ with $(m,m) \leq_m (j,\beta) \leq_m (n,n)$. In \cite[Proposition 1.3.3.2]{lau3}, we have proved (in a more general setting) that the $\mathbb{C}$-algebra homomorphism $\theta_{j,\beta} : R_{j,\beta}^{(j,\beta)^+} 
\rightarrow R_{j,\beta}^{(j,\beta)}$ induces a $\mathbb{C}$-algebra homomorphism $\ov{\theta_{j,\beta}} : B_{j,\beta}^{(j,\beta)^+} 
\rightarrow B_{j,\beta}^{(j,\beta)}$ which sends  $y_{i,\alpha}^{(j,\beta)^+}$ to $y_{\ia}^{(j,\beta)}$ ($(m,n) \leq_m 
(i,\alpha) <_m (j,\beta)$).

Let $\delta^{(j,\beta)^+}=det_q (\yia^{(j,\beta)^+})_{ \substack{i=i_1,\dots ,i_t \\ \alpha = \alpha_1, \dots \alpha_t } } $ be a quantum minor of $R^{(j,\beta)^+}$ such that $(i_l,\alpha_l) <_m (j,\beta)$ for all $l \in \gc 1,t \dc$. If  $\delta^{(j,\beta)^+}  \in L^{(j,\beta)^+}$, then $det_q (y_{\ia}^{(j,\beta)^+})_{ \substack{i=i_1,\dots ,i_t \\ \alpha = \alpha_1, \dots \alpha_t } } =0 $ in  $B_{j,\beta}^{(j,\beta)^+} $. Hence $\ov{\theta_{j,\beta}} \left(
det_q (y_{\ia}^{(j,\beta)^+})_{ \substack{i=i_1,\dots ,i_t \\ \alpha = \alpha_1, \dots \alpha_t } } \right) =0$ in 
$B_{j,\beta}^{(j,\beta)}$, that is, $det_q (y_{\ia}^{(j,\beta)})_{ \substack{i=i_1,\dots ,i_t \\ \alpha = \alpha_1, \dots \alpha_t } }  =0$ in $B_{j,\beta}^{(j,\beta)}$. In other words, $\delta^{(j,\beta)}:=det_q (\yia^{(j,\beta)})_{ \substack{i=i_1,\dots ,i_t \\ \alpha = \alpha_1, \dots \alpha_t } } \in L^{(j,\beta)}$. So we have just proved the following result.


\begin{prop}
\label{faceff2}
Let $(j,\beta) \in E$ with $(m,m) \leq_m (j,\beta) \leq_m (n,n)$ and consider a quantum minor 
$\delta^{(j,\beta)^+}=det_q (\yia^{(j,\beta)^+})_{ \substack{i=i_1,\dots ,i_t \\ \alpha = \alpha_1, \dots \alpha_t } } $ of $R^{(j,\beta)^+}$ such that $(i_l,\alpha_l) <_m (j,\beta)$ for all $l \in \gc 1,t \dc$. 
\\If $\delta^{(j,\beta)^+}  \in L^{(j,\beta)^+}$, then $\delta^{(j,\beta)}:=det_q (\yia^{(j,\beta)})_{ \substack{i=i_1,\dots ,i_t \\ \alpha = \alpha_1, \dots \alpha_t } } \in L^{(j,\beta)}$.
\end{prop}

\subsection{A criterion for quantum minors to belong to a prime ideal of $\sn$ that lies in the image of $\psi$.}
$ $

In this section, we keep the notations and conventions of the previous paragraph. In particular, $J$ still denotes an $\hc$-invariant prime ideal of $A=\mmpc$ and we set $L:= \tau \circ \iota (J)$. Moreover, we set 
$K:= \psi(J)$; this is an $\hc'_R$-invariant prime ideal of $\sn$. The aim of this section is to obtain a criterion 
  in order to know which quantum minor can belong to $K$.

\begin{prop}
\label{propcrit1}
Let $(j,\beta) \in E$ with $(j,\beta) \geq_m (m,m)$ and let $\delta^{(j,\beta)}=det_q (\yia^{(j,\beta)})_{ \substack{i=i_1,\dots ,i_t \\ \alpha = \alpha_1, \dots \alpha_t } } $ be a quantum minor of $R^{(j,\beta)}$. 
\\If $ (m,m) \leq_m (i_l,\alpha_l) <_m (j,\beta)$ for all $l \in \gc 1,t \dc$, then $\delta^{(j,\beta)} \not \in L^{(j,\beta)}$.
\end{prop}
\preuve We proceed by induction on $(j,\beta)$. If $(j,\beta)=(m,m)$, then there is nothing to prove. 
We now assume that $(m,m) \leq_m (j,\beta) \leq_m (n,n)$. Let $\delta^{(j,\beta)^+}=det_q (\yia^{(j,\beta)^+})_{ \substack{i=i_1,\dots ,i_t \\ \alpha = \alpha_1, \dots \alpha_t } } $ be a quantum minor of $R^{(j,\beta)^+}$ such that $ (m,m) \leq_m (i_l,\alpha_l) <_m (j,\beta)^+$ for all $l \in \gc 1,t \dc$. We need to show that $\delta^{(j,\beta)^+} \not \in L^{(j,\beta)^+}$. We assume that this is not the case, that is, we assume that $\delta^{(j,\beta)^+} \in L^{(j,\beta)^+}$ and we distinguish several cases.
\\$\bullet$ Assume that $i_t > m$ and $(i_t,\alpha_t) <_m (j,\beta)$. Since $\delta^{(j,\beta)^+} \in L^{(j,\beta)^+}$,  we deduce from Proposition \ref{faceff2} that 
$\delta^{(j,\beta)}:=det_q (\yia^{(j,\beta)})_{ \substack{i=i_1,\dots ,i_t \\ \alpha = \alpha_1, \dots \alpha_t } } \in L^{(j,\beta)}$. On the other hand, it follows from the induction hypothesis that $\delta^{(j,\beta)} \not \in L^{(j,\beta)}$. This is a contradiction.
\\$\bullet$ Assume that $i_t > m$ and $(i_t,\alpha_t) = (j,\beta)$. Then we have $j > m$ and so it follows from 
\cite[Remarques 3.2.10]{lau} that the $(m,p)$ deleting-derivations algorithm coincides with the standard one (see \cite{lau}). With this in mind, \cite[Proposition 4.2.1]{c2} leads to the following equality 
$$\delta^{(j,\beta)^+}=det_q (\yia^{(j,\beta)^+})_{ \substack{i=i_1,\dots ,i_t \\ \alpha = \alpha_1, \dots \alpha_t } } =det_q (\yia^{(j,\beta)})_{ \substack{i=i_1,\dots ,i_{t-1} \\ \alpha = \alpha_1, \dots \alpha_{t-1} } }
V_{j,\beta}. $$
Recall that $\Sigma_{j,\beta}$ denotes the multiplicative system of $\fract(R)$ generated by $V_{j,\beta}$. Since $\delta^{(j,\beta)^+} \in L^{(j,\beta)^+}$, we deduce from the previous equality that 
$$det_q (\yia^{(j,\beta)})_{ \substack{i=i_1,\dots ,i_{t-1} \\ \alpha = \alpha_1, \dots \alpha_{t-1} } } V_{j,\beta}
\in L^{(j,\beta)^+} \Sigma_{j,\beta}^{-1} \cap R^{(j,\beta)}.$$
Since $L^{(j,\beta)^+} \Sigma_{j,\beta}^{-1} \cap R^{(j,\beta)}=L^{(j,\beta)}$ (see Observation \ref{obsutile}), this leads to  
$$det_q (\yia^{(j,\beta)})_{ \substack{i=i_1,\dots ,i_{t-1} \\ \alpha = \alpha_1, \dots \alpha_{t-1} } } V_{j,\beta}
\in L^{(j,\beta)}.$$
Now recall that $L^{(j,\beta)} \in \hc_R$-$\spec^{**}(R^{(j,\beta)})$. In particular, $V_{j,\beta}$ does not belong to the completely prime ideal $L^{(j,\beta)}$ of $R^{(j,\beta)}$. Hence we get 
$$det_q (\yia^{(j,\beta)})_{ \substack{i=i_1,\dots ,i_{t-1} \\ \alpha = \alpha_1, \dots \alpha_{t-1} } }  \in L^{(j,\beta)}.$$
On the other hand, we deduce from the induction hypothesis that $det_q (\yia^{(j,\beta)})_{ \substack{i=i_1,\dots ,i_{t-1} \\ \alpha = \alpha_1, \dots \alpha_{t-1} } } \not \in L^{(j,\beta)}$. And so we also get a contradiction in this case.
\\$\bullet$ Assume that $i_t \leq m$ and $(i_1,\alpha_1) <_m (j,\beta)$. Since $\delta^{(j,\beta)^+} \in L^{(j,\beta)^+}$, we deduce from Proposition \ref{faceff2} that 
$\delta^{(j,\beta)}:=det_q (\yia^{(j,\beta)})_{ \substack{i=i_1,\dots ,i_t \\ \alpha = \alpha_1, \dots \alpha_t } } \in L^{(j,\beta)}$.  On the other hand, it follows from the induction hypothesis that $\delta^{(j,\beta)} \not \in L^{(j,\beta)}$. This is a contradiction.
\\$\bullet$ Assume $i_t \leq m$ and $(i_1,\alpha_1) = (j,\beta)$, then we have $j \leq m $ and so it follows from 
\cite[Remarques 3.2.10]{lau} that the $(m,p)$ deleting-derivations algorithm coincides with the inverse one (see \cite{lau}). With this in mind, \cite[Proposition 2.5.6]{lau} leads to the following equality 
$$\delta^{(j,\beta)^+}=det_q (\yia^{(j,\beta)^+})_{ \substack{i=i_1,\dots ,i_t \\ \alpha = \alpha_1, \dots \alpha_t } } =det_q (\yia^{(j,\beta)})_{ \substack{i=i_2,\dots ,i_{t} \\ \alpha = \alpha_2, \dots \alpha_{t} } }
V_{j,\beta}. $$
Since $\delta^{(j,\beta)^+} \in L^{(j,\beta)^+}$, we deduce from the previous equality that 
$$det_q (\yia^{(j,\beta)})_{ \substack{i=i_2,\dots ,i_{t} \\ \alpha = \alpha_2, \dots \alpha_{t} } } V_{j,\beta}
\in L^{(j,\beta)^+} \Sigma_{j,\beta}^{-1} \cap R^{(j,\beta)}.$$
Since $L^{(j,\beta)^+} \Sigma_{j,\beta}^{-1} \cap R^{(j,\beta)}=L^{(j,\beta)}$ (see Observation \ref{obsutile}), this leads to  
$$det_q (\yia^{(j,\beta)})_{ \substack{i=i_2,\dots ,i_{t} \\ \alpha = \alpha_2, \dots \alpha_{t} } } V_{j,\beta}
\in L^{(j,\beta)}.$$
Now recall that $L^{(j,\beta)} \in \hc_R$-$\spec^{**}(R^{(j,\beta)})$. In particular, $V_{j,\beta}$ does not belong to the completely prime ideal $L^{(j,\beta)}$ of $R^{(j,\beta)}$. Hence we get 
$$det_q (\yia^{(j,\beta)})_{ \substack{i=i_2,\dots ,i_{t} \\ \alpha = \alpha_2, \dots \alpha_{t} } }  \in L^{(j,\beta)}.$$
On the other hand, it follows from the induction hypothesis that $det_q (\yia^{(j,\beta)})_{ \substack{i=i_2,\dots ,i_{t} \\ \alpha = \alpha_2, \dots \alpha_{t} } } \not \in L^{(j,\beta)}$. And so we get a contradiction.
\\$ $

To resume, we have obtained a contradiction in all cases. This proves that $\delta^{(j,\beta)^+} \not \in L^{(j,\beta)^+}$, as desired. \fin
\\$ $

We are now able to prove the following result that provides a useful tool to prove that a quantum minor does not belong to an ideal in the image of the embedding $\psi$.

\begin{theo}
\label{critere}
Let $J  \in \hc\mbox{-} \spec(A)$ and set $K:= \psi(J) \in \hc'_R \mbox{-} \spec(\sn)$. Let $\delta=det_q (\xia)_{ \substack{i=i_1,\dots ,i_t \\ \alpha = \alpha_1, \dots \alpha_t } } $ be a quantum minor of $\sn$.
\\If $\delta \in K$, then there exists $l \in \gc 1,t \dc $ such that $(i_l, \alpha_l) \in \gc 1,m \dc \times \gc m+1,n \dc$.
\end{theo}
\preuve Let $\delta=det_q (\xia)_{ \substack{i=i_1,\dots ,i_t \\ \alpha = \alpha_1, \dots \alpha_t } } $ be a quantum minor of $\sn$ that belongs to $K$. Set $L:= \tau \circ \iota (J)$. It follows from the construction of $\psi$ that, with the notation of section \ref{sectionSn}, we have $L= \theta(\widehat{K}) \cap R$, where $\theta : \sn[z,z^{-1}] \rightarrow \gn$ denotes the isomorphism constructed by Levasseur and Stafford (see section \ref{sectionSn}). In \cite[Observation 3.4.8]{lau}, we have computed the image of a quantum minor of $\sn$ by this isomorphism $\theta$. These computations together with the fact that $\delta$ belongs to 
$K$ show that $det_q (\yia)_{ \substack{i=i_1,\dots ,i_t \\ \alpha = \alpha_1, \dots \alpha_t } } \in \theta(\widehat{K}) \cap R=L$. Next Proposition \ref{propcrit1} (with $(j,\beta)=(n,n+1)$) shows that there must exist $l \in \gc 1,t \dc $ such that $(i_l, \alpha_l) <_m (m,m)$. In other words, there exists $l \in \gc 1,t \dc $ such that $(i_l, \alpha_l) \in \gc 1,m \dc \times \gc m+1,n \dc$. \fin
\\$ $


\subsection{A transfert result.}
$ $

The aim of this paragraph is to prove that, if $J \in \hc$-$\spec(A)$ contains all quantum minors of $A$ of a given size $r$, then $\psi(J)$ contains all $r \times r$ quantum minors of the matrix $(\xia)_{(\ia) \in \gc 1,m \dc \times \gc m+1,n \dc}$.

\begin{prop}
\label{lemdebase}
Let $J \in \hc$-$\spec(A)$ and $r \in \gc 0, \min(m,p) \dc$.  Assume that $J$ contains all $r \times r$ quantum minors of $A$.
\\Then $ det_q (\xia)_{ \substack{i=i_1,\dots ,i_r \\ \alpha = \alpha_1, \dots \alpha_r } } \in  \psi(J) $  for all $ 1 \leq i_1 < \dots < i_r \leq m$ and $m+1 \leq \alpha_1 < \dots < \alpha_r \leq n$.
\end{prop}
\preuve As in the previous sections, we set $L:= \tau \circ \iota (J)$; this is an $\hc_R$-invariant prime ideal of $R$ that does not contain the quantum determinant $\Delta$. Further, still as in the previous sections, we set $L^{(m,m)}:= \iota (J)$ and, for all $(j,\beta) \in E$ with $(m-1,m) \leq_m (j,\beta) \leq_m (n,n+1)$, $L^{(j,\beta)}:= \tau_{(j,\beta)^-} \circ \dots \circ \tau_{m,m} \circ \iota (J)$ where $(j,\beta)^-$ denotes the greatest element (relative to $\leq_m$) of the set $\{(\ia) \in E \mid (\ia) < (j,\beta)\}$; $L^{(j,\beta) }$ is an $\hc_R$-invariant prime ideal of $R^{(j,\beta)}$. Recall (see \ref{sectionLjbeta}) that, for all $(j,\beta) \in E$ with 
$(m,m) \leq_m (j,\beta) \leq_m (n,n)$, we have $L^{(j,\beta)} \Sigma_{j,\beta}^{-1} = L^{(j,\beta)^+} \Sigma_{j,\beta}^{-1}$, where $\Sigma_{j,\beta} = \{ V_{j,\beta}^i \mid i \in \mathbb{N} \}$
.

First, we show that, for all $(j,\beta) \in E$ with $(j,\beta) \geq (m,m)$, $det_q(\yia^{(j,\beta)} )_{ \substack{i=i_1,\dots ,i_r \\ \alpha = \alpha_1, \dots \alpha_r } } \in  L^{(j,\beta)}   $  for all $ 1 \leq i_1 < \dots < i_r \leq m$ and $m+1 \leq \alpha_1 < \dots < \alpha_r \leq n$. We proceed by induction on $(j,\beta)$.

Assume $(j,\beta)=(m,m)$ and recall that we have identified $A$ to the subalgebra of $R^{(m,m)}$ generated by those 
$ \yia^{(m,m)}$ such that $(\ia) \in \gc 1,m \dc \times \gc m+1,n \dc$. Since $J $ contains all $r \times r$ quantum minors of $A$, we have $ det_q (\yia^{(m,m)})_{\substack{i=i_1,\dots ,i_r \\ \alpha = \alpha_1, \dots \alpha_r }   } \in  J $  for all $ 1 \leq i_1 < \dots < i_r \leq m$ and $m+1 \leq \alpha_1 < \dots < \alpha_r \leq n$. Thus, since $J \subseteq L^{(m,m)}$, we get $ det_q (\yia^{(m,m)})_{ \substack{i=i_1,\dots ,i_r \\ \alpha = \alpha_1, \dots \alpha_r }  } \in  L^{(m,m)} $  for all $ 1 \leq i_1 < \dots < i_r \leq m$ and $m+1 \leq \alpha_1 < \dots < \alpha_r \leq n$, as required.

Assume now that $(j,\beta) \in E$ with $(j,\beta) \neq (n,n+1)$. Let $ \delta=det_q (\yia^{(j,\beta)^+})_{ \substack{i=i_1,\dots ,i_r \\ \alpha = \alpha_1, \dots \alpha_r }  } $ be a quantum minor of $R^{(j,\beta)^+}$ with $ 1 \leq i_1 < \dots < i_r \leq m$ and $m+1 \leq \alpha_1 < \dots < \alpha_r \leq n$. In order to prove that $\delta$ belongs to $L^{(j,\beta)^+}$, we distinguish two cases. 
\\$\bullet$ If $j > m$, it follows from \cite[Remarques 3.2.10]{lau} that the $(m,p)$ deleting-derivations algorithm coincides with the standard one (see \cite{lau}). Hence we deduce from \cite[Proposition 2.2.8]{lau} that $\delta$ is 
a right linear combination with coefficients in $R^{(j,\beta)} \Sigma_{j,\beta}^{-1}$ of quantum minors $ det_q (\yia^{(j,\beta)})_{ \substack{i=i_1,\dots ,i_r \\ \alpha = \beta_1, \dots \beta_r } } $ of $R^{(j,\beta)}$ with $ \beta_k \geq \alpha_k$ for all $k \in \gc 1,r \dc$. Now it follows from the induction hypothesis that all these quantum minors belong to $L^{(j,\beta)}$, so that $\delta \in L^{(j,\beta)} \Sigma_{j,\beta}^{-1}=L^{(j,\beta)^+} \Sigma_{j,\beta}^{-1}$. Hence $\delta \in L^{(j,\beta)^+} \Sigma_{j,\beta}^{-1} \cap R^{(j,\beta)^+} =L^{(j,\beta)^+}$, as desired.
\\$\bullet$ If $j \leq m$, it follows from \cite[Remarques 3.2.10]{lau} that the $(m,p)$ deleting-derivations algorithm coincides with the inverse one (see \cite{lau}). Note that, in \cite{lau}, we have not established an analogue of \cite[Proposition 2.2.8]{lau} for the inverse deleting-derivations algorithm. However, using  \cite[Proposition 2.2.8]{lau} and \cite[Corollaire 2.4.3]{lau}, it is easy to obtain that, in this case, $\delta$ is 
a right linear combination with coefficients in $R^{(j,\beta)} \Sigma_{j,\beta}^{-1}$ of quantum minors $ det_q (\yia^{(j,\beta)})_{ \substack{i=j_1,\dots ,j_r \\ \alpha = \alpha_1, \dots \alpha_t } } $ of $R^{(j,\beta)}$ with $ j_k \leq i_k$ for all $k \in \gc 1,r \dc$. Now, following the route sketched in the previous case, we also get $ \delta \in L^{(j,\beta)^+}$.

To resume, we have just proved that, for all $(j,\beta) \in E$ with $(j,\beta) \geq (m,m)$, every quantum minor $ det_q (\yia^{(j,\beta)})_{ \substack{i=i_1,\dots ,i_r \\ \alpha = \alpha_1, \dots \alpha_r } } $ with $1 \leq i_1 < \dots < i_r \leq m $ and $m+1 \leq \alpha_1 < \dots < \alpha_r \leq n $ belongs to $  L^{(j,\beta)}  $. In particular, for $(j,\beta)=(n,n+1)$, we obtain that every quantum minor $ det_q (\yia)_{ \substack{i=i_1,\dots ,i_r \\ \alpha = \alpha_1, \dots \alpha_r } } $ with $1 \leq i_1 < \dots < i_l \leq m $ and $m+1 \leq \alpha_1 < \dots < \alpha_l \leq n $ belongs to $  L$.

Finally, let $\delta=det_q (\xia)_{ \substack{i=i_1,\dots ,i_r \\ \alpha = \alpha_1, \dots \alpha_r } } $ be a quantum minor of $\sn$ with $1 \leq i_1 < \dots < i_r \leq m $ and $m+1 \leq \alpha_1 < \dots < \alpha_r \leq n $. It follows from the construction of $\psi$ that, with the notation of section \ref{sectionSn}, we have $\psi(J)=  \theta^{-1}(L_{\Delta}) \cap \sn$, where $\theta : \sn [z,z^{-1}] \rightarrow \gn$ denotes the isomorphism constructed by Levasseur and Stafford (see section \ref{sectionSn}). In \cite[Observation 3.4.8]{lau}, we have computed the image of a quantum minor of $\sn$ by this isomorphism $\theta$. Since $ det_q (\yia)_{ \substack{i=i_1,\dots ,i_r \\ \alpha = \alpha_1, \dots \alpha_r } } $ belongs to $L$, we deduce from these computations that $\delta$ belongs to the completely prime ideal $\psi(J)$, as desired. \fin


\section{An ordered bijection between $ \hc$-$\spec(A)$ and $\mathcal{S}$.}
$ $

In section \ref{sectionpsi}, we have constructed an order-preserving embedding $\psi : \hc \mbox{-} \spec (A) \rightarrow \hc'_R \mbox{-} \spec (\sn )$. Our first aim in this section is to determine the image of $\psi$. Recall that Brown and Goodearl have determined the set $\hc'_R \mbox{-} \spec (\sn )$. More precisely, they have shown (see Proposition \ref{propHprimesln}) that $\hc'_R \mbox{-} \spec (\sn ) = \{ I_w \mid w \in S_n \times S_n \}$, where each $I_w$ is generated by some quantum minors. In this section, we prove that the image of the embedding $\psi$ is exactly the set of those $I_{w_0,w_0\sigma}$ such that $\sigma \in  \mathcal{S}$  (where $w_0$ still denotes the longest element of $S_n$ and $\mathcal{S}=\{ \sigma' \in S_n \mid -p \leq i-\sigma'(i) \leq m \mbox{ for all } i \in \gc 1, n \dc \}$ is the sub-poset of the (reverse) Bruhat order of $S_n$ that we have studied in section \ref{sectionbruhatS}). To do this, we use the criteria obtained in the previous section. These criteria allow us to prove that the image of $\psi$ is contained in the set of those $I_{w_0,w_0\sigma}$ such that $\sigma \in  \mathcal{S}$. Next, since we have already proved (see Corollary \ref{cardok}) that the sets $\hc$-$\spec(A)$ and $\mathcal{S}$ have the same cardinality, we obtain that the image of $\psi$ is actually equal to the set of those $I_{w_0,w_0\sigma}$ such that $\sigma \in  \mathcal{S}$. As a consequence, $\psi$ induces an ordered bijection between $\hc$-$\spec(A)$ and the set of those $I_{w_0,w_0\sigma}$ such that $\sigma \in  \mathcal{S}$. Next we prove that the inverse bijection is also an order-preserving one whose construction is algorithmic. Naturally we derive from this bijection  an ordered bijection $\Xi$ from the sub-poset $\mathcal{S}$ of the (reverse) Bruhat order of $S_n$ and $\hc$-$\spec(A)$, whose construction is also algorithmic. The existence of such a bijection has been conjectured by Goodearl and Lenagan.

In the last section, we investigate, in the case where $m=p$, the combinatoric of the set of rank $t$ $\hc$-invariant prime ideals of $A$, that is, the combinatoric of the set of those $\hc$-invariant prime ideals of $A$ that contain all $(t+1) \times (t+1)$ quantum minors, but not all $t \times t$ quantum minors. More precisely, we prove the following result. Imagine that there is a barrier between positions $m$ and $m+1$. Then rank $t$ $\hc$-invariant prime ideals of $A$ are those $\hc$-invariant prime ideals of $A$
that correspond via the bijection $\Xi^{-1}$ to  $n$-permutations $\sigma \in \mathcal{S}$ such that the number of integers that are moved by $\sigma$ from the right to the left of this barrier is exactly $m-t$. This result was also conjectured by Goodearl and Lenagan.

\subsection{The image of the embedding $\psi$.}
\label{sectionimagepsi}
$ $

 Recall (see Proposition \ref{propHprimesln}) that  $\hc'_R \mbox{-}\spec(\sn) = \{I_w \mid w=(w_+,w_-) \in S_n \times S_n\}$, where $I_w$ denotes the ideal of $\sn $ generated by the quantum minors $c_{j,y}^{+}:=det_q(\xia)_{(\ia) \in y(\gc 1, j \dc ) \times \gc 1,j \dc}$ with $j\in \gc 1,n-1 \dc $, $y \in S_n$ and $y \nleq_j w_{+}$, together with the quantum minors $c_{j,y}^-:=det_q(\xia)_{(\ia) \in y(\gc j+1,n \dc ) \times \gc j+1,n \dc}$ with $j\in \gc 1,n-1 \dc $, $y \in S_n$ and $y \nleq_j w_{-}$ (see Notation \ref{notationsideauxhodgeslevasseur}).

We first use the criteria obtained in the previous section in order to show that  
$\psi ( \hc \mbox{-} \spec(A)) \subseteq \{ I_{w_0,w_0\sigma} \mid \sigma \in \mathcal{S} \}$ (where $w_0$ still denotes the longest element of $S_n$ and $\mathcal{S}$ the sub-poset of the (reverse) Bruhat order of $S_n$ that we have studied in section \ref{sectionbruhatS}). 
\\$ $

Throughout this section, we will adopt the following notations. If $y \in S_n$ and $j \in \gc 1,n-1 \dc$, we set $y(\gc 1,j \dc)=\{y_1^j < \dots < y_j^j\}$. Clearly, we have $y_l^j \geq l$ for all $l \in \gc 1,j \dc$.\\

\begin{lem}
\label{lemjw1}
Let $w=(w_+,w_-) \in S_n \times S_n$. If $I_w \in \psi ( \hc \mbox{-} \spec(A))$, then $w_+=w_0$.
\end{lem}
\preuve Assume that this is not 
the case, that is, assume that $w_+ \neq w_0$. Then $ w_+ < w_0$ and it follows from Proposition \ref{propbruhat} that, there exists $j \in \gc 1,n-1 \dc$ such that $w_+ <_j w_0$. Hence $w_0 \nleq_j w_+$, and so we have 
$$ c_{j,w_0}^+:=det_q(\xia)_{(\ia) \in w_0(\gc 1, j \dc ) \times \gc 1,j \dc} \in I_w,$$
that is, 
$$ c_{j,w_0}^+:=det_q(\xia)_{(\ia) \in \gc n+1-j, n \dc  \times \gc 1,j \dc} \in I_w,$$
Next, since $I_w \in \psi ( \hc \mbox{-} \spec(A))$, we deduce from Theorem \ref{critere} that there must exist $i \in \gc 1,j \dc$ such that $(n-j+i,i)  \in \gc 1,m \dc \times \gc m+1,n \dc$. This is a contradiction. \fin

\begin{lem}
\label{lemjw2}
Let $w=(w_+,w_-) \in S_n \times S_n$. If $I_w \in \psi ( \hc \mbox{-} \spec(A))$, then $w_-(n-t) \leq m+1+t$ for all $t \in \gc 0,p-2 \dc$.
\end{lem}
\preuve We proceed by induction on $t$. 
First, assume that $w_-(n-0) \nleq m+1+0$, that is, $w_-(n) \geq m+2$. Set 
$w_-(\gc 1,n-1 \dc) = \{w_1 < \dots < w_{n-1} \}$. Since $w_-(n) \geq m+2$, we have $w_{m+1}=m+1$.  
Now we consider the transposition $y:=s_{m+1,n} \in S_n$. 
It is clear that $y(\gc 1,n -1 \dc) = \{1,\dots,m,m+2,\dots,n \}$. Hence $y_{m+1}^{n-1}=m+2 > w_{m+1}=m+1$ and so $y \nleq_{n-1} w_-$. Thus we have 
$c_{n-1,y}^-=X_{m+1,n} \in I_w$. However, since $I_w \in \psi ( \hc \mbox{-} \spec(A))$, we deduce from Theorem \ref{critere} that $X_{m+1,n}$ can not belong to $I_w$. So we get a contradiction and thus we have proved that 
$w_-(n-0) \leq m+1+0$.

Let now $t \in \gc 1,p-2 \dc$ and assume that $w_-(n-t) \nleq m+1+t$, that is, $w_-(n-t) \geq m+2+t$. 
Set $j:= n-1-t \in \gc m+1,n-2 \dc$ and $w_-(\gc 1,j \dc) = w_-(\gc 1,n-1-t \dc) := \{w_1 < \dots < w_{n-1-t} \}$. It follows from the induction hypothesis that $w_-(n) \leq m+1$,..., $w_-(n-t+1) \leq m+t$. 
Thus, since  $w_-(n-t) \geq m+2+t$, we get
\begin{eqnarray*}
w_-(\gc 1,j \dc) & = & \left( \gc 1,m+t \dc \setminus \left\{w_-(i) \mid i \in \gc n+1-t , n \dc \right\} \right) \\
 & & \bigcup \left( \gc m+t+1,n \dc \setminus \{w_-(n-t)\} \right) .
\end{eqnarray*}
In particular, we have $w_{m+1}=m+t+1$. Now, let $y$ be the permutation of $S_n$ defined by 
$y(k)=k$ if $k \leq m$ and $y(k)=m+n+1-k $ if $k \geq m+1$. Clearly, we have 
$y(\gc 1,j \dc)=y(\gc 1,n-1-t \dc)=\gc 1,m \dc \cup \gc m+t+2, n \dc$. Hence $y_{m+1}^j=m+t+2 > m+t+1=w_{m+1}$ and so $y \nleq_{j} w_-$. Thus we have 
$c_{j,y}^- = det_q(\xia)_{(\ia) \in \gc m+1,m+1+t \dc \times \gc n-t,n \dc} \in I_w$. Next, since $I_w \in \psi ( \hc \mbox{-} \spec(A))$, we deduce from Theorem \ref{critere} that there must exist $l \in \gc 0,t \dc$ such that 
$(m+1+l, n-t+l) \in \gc 1,m \dc \times \gc m+1,n\dc$. This is a contradiction and thus we have proved that 
$w_-(n-t) \leq m+1+t$. \fin

\begin{lem}
\label{lemjw3}
Let $w=(w_+,w_-) \in S_n \times S_n$. If $I_w \in \psi ( \hc \mbox{-} \spec(A))$, then $w_-(t) \geq m+1-t$ for all $t \in \gc 1,m-1 \dc$.
\end{lem}
\preuve We proceed by induction on $t$. 
First, assume that $w_-(1) \ngeq m+1-1$, that is, assume that $w_-(1) \leq m-1$. Set $y:=s_{1, m} \in S_n$. Clearly, we have $y \nleq_{1} w_-$, so that 
$c_{1,y}^-=det_q(\xia)_{(\ia) \in \{1,\dots,m-1,m+1,\dots,n \} \times \gc 2,n \dc} \in I_w$. Next, since $I_w \in \psi ( \hc \mbox{-} \spec(A))$, we deduce from Theorem \ref{critere} that either there exists $l \in \gc 1,m-1 \dc$ such that $(l,l+1) \in \gc 1,m \dc \times \gc m+1,n\dc$ or there exists $l \in \gc m+1,n \dc$ such that $(l,l)\in \gc 1,m \dc \times \gc m+1,n\dc$. In both cases, we obtain a contradiction and thus we have proved that 
$w_-(1) \geq m$.

Let now $t \in \gc 2,m-1 \dc$ and assume that $w_-(t) \ngeq m+1-t$, that is, $w_-(t) \leq m-t$. 
Set $w_-(\gc 1,t \dc) =  \{w_1 < \dots < w_{t} \}$. It follows from the induction hypothesis that $w_-(1) \geq m$,..., $w_-(t-1) \geq m+2-t$. 
Thus, since  $w_-(t) \leq m-t$, we get $w_{1}=w_-(t) \leq m-t$. Now, let $y$ be the permutation of $S_n$ defined by 
$y(k)=m+1-k$ if $k \leq m$ and $y(k)=k $ if $k \geq m+1$. Clearly, we have 
$y(\gc 1,t \dc)= \gc  m+1-t,m \dc$. Hence $y_{1}^t= m+1-t > m-t \geq w_{1}$ and so $y \nleq_{t} w_-$. Thus $c_{t,y}^- = det_q(\xia)_{(\ia) \in \{1,\dots, m-t,m+1,\dots,n \} \times \gc t+1,n \dc}$ belongs to $I_w$. However, since $I_w \in \psi ( \hc \mbox{-} \spec(A))$, we deduce from Theorem \ref{critere} that either there exists $l \in \gc 1,m-t \dc$ such that 
$(l, t+l) \in \gc 1,m \dc \times \gc m+1,n\dc$ or there exists $l \in \gc m+1,n \dc$ such that $(l,l)\in \gc 1,m \dc \times \gc m+1,n\dc$. In both cases, we obtain a contradiction and thus we have proved that 
$w_-(t) \geq m+1-t$. \fin
\\$ $

Recall that the set $\mathcal{S}=\{ \sigma \in S_n \mid -p \leq i - \sigma(i) \leq m \mbox{ for all } i \in \gc 1,n \dc \}$ is a sub-poset of the Bruhat order of $S_n$ whose cardinality is equal to the cardinality of $\hc$-$\spec(A)$ (see Proposition \ref{descriptionS} and Corollary \ref{cardok}) .

\begin{theo}
\label{thetheo}
Let $w=(w_+,w_-) \in S_n \times S_n$.
\\$I_w \in \psi ( \hc \mbox{-} \spec(A))$ if and only if $w_+=w_0$ and $w_0 w_- \in \mathcal{S}$.
\end{theo}
\preuve First we easily deduce from Lemmas \ref{lemjw1}, \ref{lemjw2} and \ref{lemjw3} that, if $I_w \in \psi ( \hc \mbox{-} \spec(A))$, then $w_+=w_0$ and $w_0 w_- \in \mathcal{S}$. Hence 
$\psi ( \hc \mbox{-} \spec(A)) \subseteq \{I_w \mid w_+=w_0 \mbox{ and } w_0 w_- \in \mathcal{S} \}$. On the other hand, since $\psi $ is an embedding, we have 
$$ \mid \psi ( \hc \mbox{-} \spec(A)) \mid = \mid  \hc \mbox{-} \spec(A) \mid$$
and so we deduce from Proposition \ref{cardinalityS} that 
$$ \mid \psi ( \hc \mbox{-} \spec(A)) \mid = \mid \mathcal{S} \mid = \mid \{I_w \mid w_+=w_0 \mbox{ and } w_0 w_- \in \mathcal{S} \} \mid .$$
Thus $\psi ( \hc \mbox{-} \spec(A)) \subseteq \{I_w \mid w_+=w_0 \mbox{ and } w_0 w_- \in \mathcal{S} \}$ and these two sets have the same cardinality. This implies that 
$\psi ( \hc \mbox{-} \spec(A)) = \{I_w \mid w_+=w_0 \mbox{ and } w_0 w_- \in \mathcal{S} \}$, as desired. \fin
\\$ $

Since $\psi : \hc \mbox{-} \spec(A) \rightarrow \hc'_R \mbox{-} \spec(\sn)$ is an ordered embedding, $\psi$ induces an ordered bijection $\widetilde{\psi}$ from  $\hc \mbox{-} \spec (A)$ onto $\psi \left( \hc \mbox{-} \spec (A) \right)$. Recall (see section \ref{sectionpsi}) that the inverse of $\widetilde{\psi}$ is also an ordered bijection and that, with the notation of section \ref{sectionpsi}, we have 
$$\widetilde{\psi}^{-1}(K) =  \psi_{(m,m)}\left( \varphi(K) \cap R \right) \cap A.$$ 
An immediate corollary of Theorem \ref{thetheo} is the following result.

\begin{cor}
The map $\widetilde{\psi} : \hc \mbox{-} \spec(A) \rightarrow \{I_w \mid w_+=w_0 \mbox{ and } w_0 w_- \in \mathcal{S} \}$ is an order-preserving bijection whose inverse is also ordered.
\end{cor}

Now, recall (see Proposition \ref{bijsn}) that the map $\chi : S_n \times S_n \rightarrow \{I_w \mid w \in S_n \times S_n \}$ defined by $\chi(w)=I_w$ is an order-reversing bijection. Moreover, it is well known that the map $ f : S_n \rightarrow S_n$ defined by 
$f(\sigma)=w_0 \sigma$ is an order-reversing involution (relative to the (reverse) Bruhat order on $S_n$). Now, composing the three map $\widetilde{\psi}^{-1}$, $ \chi$ and $f$ leads to our main result.

\begin{theo}
\label{thetheo2}
The map $\Xi : \mathcal{S} \rightarrow  \hc \mbox{-} \spec(A)$ defined by 
$\Xi ( \sigma ) = \widetilde{\psi}^{-1}(I_{w_0, w_0 \sigma})$ is an ordered bijection whose inverse is also 
ordered.
\end{theo}

Note that, in the case where $p=m$, $\mathcal{S}$ is the sub-poset of the Bruhat ordering of $S_{2m}$ consisting of the permutations that move any integer by no more than $m$ positions. Thus Theorem \ref{thetheo2} furnishes an ordered bijection from $\hc$-$\spec(A)$ into this sub-poset, whose construction is actually algorithmic. The existence of such an ordered bijection has been conjectured by Goodearl and Lenagan.

Very recently, Brown, Goodearl and Yakimov (see \cite{bgy}) have investigated the geometry of the Poisson variety $\mathcal{M}_{m,p}(\mathbb{C})$, where the Poisson structure is the standard one that can be obtained from commutators of $\mmpc$. They prove in particular that  there exists a bijection between the sub-poset $\mathcal{S}$ and the set of $\hc$-orbits of symplectic leaves in $\mathcal{M}_{m,p}(\mathbb{C})$. Hence we deduce from Theorem \ref{thetheo2} the existence of a "natural" bijection between $\hc$-$\spec(\mmpc)$ and the set of $\hc$-orbits of symplectic leaves in $\mathcal{M}_{m,p}(\mathbb{C})$. Naturally this indicates that it should exist a bijection between the primitive ideals of $\mmpc$ and the symplectic leaves of $\mathcal{M}_{m,p}(\mathbb{C})$. Recall that, in the case of $SL_n(\mathbb{C})$, such kind of bijection was constructed by Hodges and Levasseur (see \cite{hl1,hl2}) whose results were extended to an arbitrary connected, simply connected, semi-simple complex algebraic group by Joseph (see \cite{jos}).

\subsection{The rank $t$ $\hc$-primes of $\mmc$.}
\label{sectionrangtperm}
$ $

Throughout this section, we assume that $m=p$, so that $A=\mmc$. The aim of this section is to recognize which permutations of 
$S_{2m}$ correspond via the bijection $\Xi$ to rank $t$ $\hc$-invariant prime ideals of $A$, that is, to $\hc$-invariant prime ideals of $A$ that contain all $(t+1) \times (t+1)$ quantum minors, but not all $t \times t$ quantum minors. In fact, we will establish the following result. Let $\sigma$ be a $2m$-permutation belonging to $\mathcal{S}$. Imagine that there is a barrier between positions $m$ and $m+1$. Then $\Xi (\sigma)$ is an $\hc$-invariant prime ideal of $A$ that has rank $m-t$ if and only if the number of integers that are moved by $\sigma$ from the right to the left of this barrier is exactly $t$. To prove this result, we will proceed in two steps. 
First, we prove that the number of permutations $\sigma \in \mathcal{S}$ such that the number of integers that are moved by $\sigma$ from the right to the left of this barrier is exactly $t$ is the same as the number of $\hc$-invariant prime ideals of $A$ that have rank $m-t$. Next we will show that, if $\Xi(\sigma)$ has rank $m-t$, then the number of integers that are moved by $\sigma$ from the right to the left of this barrier is exactly $t$.
\\$ $

\begin{defi} Let $r \in \gc 0,m \dc$. An $\hc$-invariant prime ideal $J$ of $A$ has rank $r$ 
if $J$ contains all $(r+1) \times (r+1)$ quantum minors but not all $r \times r$ quantum minors.
\\As in \cite[3.6]{glen1}, we denote by $\hc$-$\spec^{[r]} (A)$ 
the set of rank $r$ $\hc$-invariant prime ideals of $A$.
\end{defi}
$ $

Note that there is only one element in $\hc$-$\spec^{[0]} (A)$: $\langle \yia \mbox{ $\mid$ } 
(\ia) \in \gc 1,m \dc^2 \rangle$, the augmentation ideal of $A$. Further, 
Goodearl and Lenagan have observed (see \cite{glen1}, 3.6) that 
$\mid \hc \mbox{-}\spec^{[1]} (A) \mid \ =(2^m -1)^2$ and $\mid \hc \mbox{-}\spec^{[m]} (A) \mid \ =(m!)^2$. 
Next we have given a formula for the number of rank $r$ $\hc$-invariant prime ideals of $A$ (see \cite[Theorem 3.11]{lau4}).

\begin{theo}
\label{numberrankt} 
If $r \in \gc 0,m \dc$, then we have $\mid \hc\mbox{-}\spec^{[r]} (A) \mid \ = 
\left( r! S(m+1,r+1) \right)^2$, where $S(m+1,r+1)$ still denotes the Stirling number of second kind associated to 
$m+1$ and $r+1$.
\end{theo}

Recall that $m=p$, $\mathcal{S}= \{ \sigma \in S_{2m} \mid -m \leq  \sigma(i) -i  \leq m \mbox{ for all } i \in \gc 1,2m \dc \}$, and imagine that there is a barrier between positions $m$ and $m+1$. Recall (see \ref{sectionST}) that $\mathcal{S}_t$ denotes the set of those $2m$-permutations $\sigma \in \mathcal{S}$ such that the number of integers that are moved by $\sigma$ from the right to the left of this barrier is exactly $t$ ($t \in \gc 0,m \dc$). Since $\left[ (m-t)! S(m+1,m-t+1) \right]^2$ is the number of rank $(m-t)$ $\hc$-invariant prime ideals of $A$ 
(see Theorem \ref{numberrankt}) and since this is also the cardinality of $\mathcal{S}_t$ (see Proposition \ref{cardinalSt}), we get the following result.

\begin{cor}
\label{coregalite}
For all $t \in \gc 0,m \dc$, we have $\mid \mathcal{S}_t \mid = \mid \hc\mbox{-}\spec^{[m-t]} (A) \mid $.
\end{cor}
$ $

Before proving that $\Xi(\mathcal{S}_t)=\hc\mbox{-}\spec^{[m-t]} (A)$, we first establish that every rank $(m-t)$ $\hc$-prime of $A$ belongs to $\mathcal{S}_r$ for a certain $r$ greater than or equal to $t$.

\begin{lem}
\label{lemma2}
Let $\sigma \in \mathcal{S}$ and assume that $\Xi(\sigma)$ has rank $m-t$ ($t \in \gc 1,m \dc$).
\\Then there exist $m+1 \leq j_1 < \dots < j_t \leq 2m$ such that $\sigma (j_k ) \leq m$ for all $k \in \gc 1,t \dc$.    
\end{lem}
\preuve For convenience, if $I$ and $\Gamma$ are two non-empty subsets of $\gc 1, 2m \dc$ with $\mid I \mid =\mid \Gamma \mid$, then we set $[I \mid \Gamma ] : = det_q (\xia)_{ (\ia)  \in I \times \Gamma }$.

Naturally, it is sufficient to show that, for all $l \in \gc 1,t \dc$, there exist $j_1,\dots,j_l \in \gc m+1,2m \dc$ such that $j_u \neq j_v$ if $u \neq v$ and $\sigma(j_u) \leq m$ for all $u \in \gc 1,l \dc$. We proceed by induction on $l$.

Assume that $l=1$. Since $\Xi(\sigma)$ has rank $m-t$, it follows from Proposition \ref{lemdebase} that 
$[\{1,\dots ,m\} \mid \{ m+1, \dots 2m \}]$ belongs to $I_{w_0,w_0\sigma}$. On the other hand, it follows from \cite[Theorem 2.2.3]{hl1} that $[ w_0\sigma ( \gc m+1,2m \dc ) \mid \gc m+1,2m \dc ]$ does not belong to $I_{w_0,w_0\sigma}$. Hence $w_0\sigma ( \gc m+1,2m \dc ) \neq \gc 1,m \dc$. In other words, 
there exists $j_1 \in \gc m+1,2m \dc$ such that $w_0\sigma (j_1) \geq m+1$, that is, $\sigma (j_1) \leq m$, as required.

Let now $l \in \gc 1,t-1 \dc$ and assume that there exist $j_1,\dots,j_l \in \gc m+1,2m \dc$ such that $j_u \neq j_v$ if $u \neq v$ and $\sigma(j_u) \leq m$ for all $u \in \gc 1,l \dc$. If, for all $\Gamma \subseteq \gc m+1, 2m \dc$ with $\mid \Gamma \mid = l$,  $[ w_0\sigma ( \gc m+1,2m \dc \setminus\{j_1,\dots,j_l\}) \mid \gc m+1,2m \dc \setminus \Gamma ]$ belongs to $I_{w_0,w_0\sigma}$, then, because of the $q$-Laplace relations (see, for instance, \cite[Lemma A.4]{glen1}), we should have $[ w_0\sigma ( \gc m+1,2m \dc ) \mid \gc m+1,2m \dc ] \in I_{w_0,w_0\sigma}$. However, it follows from \cite[Theorem 2.2.3]{hl1} that $[ w_0\sigma ( \gc m+1,2m \dc ) \mid \gc m+1,2m \dc ]$ does not belong to $I_{w_0,w_0\sigma}$ and so there exists $\Gamma  \subseteq \gc m+1, 2m \dc$ with $\mid \Gamma \mid = l$ such that $[ w_0\sigma ( \gc m+1,2m \dc \setminus\{j_1,\dots,j_l\}) \mid \gc m+1,2m \dc \setminus \Gamma ]$ does not belong to $I_{w_0,w_0\sigma}$. Since $m-l \geq m-t+1$ and since $\Xi(\sigma)$ has rank $m-t$, we deduce from 
Proposition \ref{lemdebase} that $I_{w_0,w_0\sigma}$ contains every $(m-l) \times (m-l)$ quantum minor $[\{i_1<\dots < i_{m-l}  \} \mid \gc m+1,2m \dc \setminus \Gamma ]$ with $i_{m-l} \leq m$. Hence we must have $w_0\sigma ( \gc m+1,2m \dc \setminus\{j_1,\dots,j_l\}) \nsubseteq \gc 1,m \dc$. In other words, there exists $j_{l+1} \in \gc m+1,2m \dc \setminus\{j_1,\dots,j_l\} $ such that $w_0\sigma (j_{l+1}) \geq m+1$, that is, $\sigma (j_{l+1}) \leq m$, as desired. \fin
\\$ $

We are now able to prove that the rank $(m-t)$ $\hc$-invariant ideals of $\mmc$ are exactly the $\hc$-invariant prime ideals of $\mmc$ that correspond via the bijection $\Xi$ to permutations of $S_{2m}$ that belong to $\mathcal{S}_t$, as conjectured by Goodearl and Lenagan.

\begin{theo}
\label{theoranktperm}
For all $t \in \gc 0,m \dc$, we have $\hc$-$\spec^{[m-t]} (\mmc) = \{ \Xi (\sigma) \mid \sigma \in \mathcal{S}_t \}$. 
\end{theo}
\preuve First, since the families $(\hc$-$\spec^{[m-t]} (\mmc))_{t \in \gc 0,m \dc}$ and 
$(\{ \Xi (\sigma) \mid \sigma \in \mathcal{S}_t \})_{t \in \gc 0,m \dc}$ form partitions of the set $\hc$-$\spec (\mnc)$, we just need to establish that 
$$\hc \mbox{-}\spec^{[m-t]} (\mmc) = \{ \Xi (\sigma) \mid \sigma \in \mathcal{S}_t \} \mbox{ for all }t \in \gc 1,m \dc.$$  We do this with the help of a decreasing induction.

Assume that $t=m$ and let $J \in \hc$-$\spec^{[0]} (\mmc)$. Then $J$ contains all $1 \times 1$ quantum minors, so that $J$ is the augmentation ideal of $A$. Since $\Xi$ is a bijection from $\mathcal{S}$ onto 
$\hc$-$\spec (\mmc)$, there exists $\sigma \in \mathcal{S}$ such that $\Xi (\sigma) =J$. Recall (see section \ref{sectionbruhatS}) that $\sigma_0$ denotes the $2m$-permutation defined by 
$$ \sigma_0 (i) = \left\{ \begin{array}{ll} m+i & \mbox{ if } i \in \gc 1,m \dc \\
i-m & \mbox{ else} 
\end{array} \right. $$
and that $\mathcal{S}= \{ \sigma' \in S_{2m} \mid \sigma' \leq \sigma_0 \}$. Since $\sigma \in \mathcal{S}$, we have 
$\sigma \leq \sigma_0$. Assume that $\sigma < \sigma_0$. Then, since $\Xi$ is an ordered bijection, we have $J=\Xi(\sigma) \subsetneq \Xi(\sigma_0)$ and the augmentation ideal $J$ is not maximal. This is a contradiction and so $\sigma=\sigma_0$. Since it is clear that $\sigma_0 \in \mathcal{S}_m$, we obtain that $J$ belongs to $  \{ \Xi (\sigma) \mid \sigma \in \mathcal{S}_m \}$, so that $\hc$-$\spec^{[0]} (\mmc) \subseteq  \{ \Xi (\sigma) \mid \sigma \in \mathcal{S}_m \}$. Now recall that we have already proved 
(see Corollary \ref{coregalite}) that these two sets have the same cardinality. Hence, they are actually equal, that 
is, $\hc$-$\spec^{[0]} (\mmc) = \{ \Xi (\sigma) \mid \sigma \in \mathcal{S}_m \}$ as desired.

Assume now that $t \in \gc 1,m-1 \dc$ and let $J \in \hc$-$\spec^{[m-t]} (\mmc)$. Since $\Xi$ is a bijection from $\mathcal{S}$ onto 
$\hc$-$\spec (\mmc)$, there exists $\sigma \in \mathcal{S}$ such that $\Xi (\sigma) =J$. Now, since $J \in \hc$-$\spec^{[m-t]} (\mmc)$, $\Xi (\sigma) $ has rank $m-t$ and so it follows from Lemma \ref{lemma2} that there exist $m+1 \leq j_1 < \dots < j_t \leq 2m$ such that $\sigma (j_k ) \leq m$ for all $k \in \gc 1,t \dc$. 
Hence $\displaystyle{\sigma \in \bigsqcup_{r=t}^m \mathcal{S}_r}$. Further, because of the induction hypothesis, $\sigma$ can not belong to $\displaystyle{\bigsqcup_{r=t+1}^m \mathcal{S}_r}$ (else the rank of $J$ will be different from $m-t$). Thus 
$\sigma \in  \mathcal{S}_t$. So we just proved that $\hc$-$\spec^{[m-t]} (\mmc) \subseteq \{ \Xi (\sigma) \mid \sigma \in \mathcal{S}_t \}$. Now recall that we have already proved (see Corollary \ref{coregalite}) that these two sets have the same cardinality. Hence, they are actually equal, that is, $\hc$-$\spec^{[m-t]} (\mmc) = \{ \Xi (\sigma) \mid \sigma \in \mathcal{S}_t \}$ as desired. \fin
\\$ $

Note that the existence of an ordered bijection from  $\mathcal{S}$ onto $\hc$-$\spec(A)$ such that Theorem \ref{theoranktperm} holds was conjectured by Goodearl and Lenagan. 
\\$ $
\begin{flushleft}
\textbf{Acknowledgments.} We wish to thank T.H. Lenagan for very helpful conversations and for bringing \cite{ves} to our attention.
\\$ $ 
\end{flushleft}



\providecommand{\bysame}{\leavevmode\hbox to3em{\hrulefill}\thinspace}
\providecommand{\MR}{\relax\ifhmode\unskip\space\fi MR }
\providecommand{\MRhref}[2]{%
  \href{http://www.ams.org/mathscinet-getitem?mr=#1}{#2}
}
\providecommand{\href}[2]{#2}

\end{document}